\documentclass[11pt]{amsart}

\usepackage[dvips]{color}

\usepackage{amssymb}
\usepackage{amsmath}
\usepackage{color}
\usepackage{epsfig}

\newtheorem{theorem}{Theorem}[section] 
\newtheorem{lemma}[theorem]{Lemma}     
\newtheorem{corollary}[theorem]{Corollary}

\newtheorem{definition}{Definition}
\newtheorem{remark}{Remark}

\newcommand{\Pp}{  {\mathbb{P}}  }

\newcommand{\R}{  {\mathbb{R}}  }
\newcommand{\C}{  {\mathbb{C}}  }

\newcommand{\al}{{\alpha}}

\author{Orlando Neto} 

\address{Universidade de Lisboa\\
 Centro de Matem\'atica e Aplica\c c\~oes Fundamentais, Av. Gama Pinto, 2 1699-003 Lisboa, Portugal \\
and Faculdade de Ci\^{e}ncias da Universidade de Lisboa, Campo Grande 1749-016, Lisboa, Portugal}
\email{orlando60@gmail.com}

\author{Pedro C. Silva}
\address{Universidade de Lisboa \\  Instituto Superior de Agronomia\\ 
 Centro de Estudos Florestais,  Tapada da Ajuda, 1349-017, Lisboa, Portugal}
\email{pcsilva@isa.utl.pt}

\title{Rigid local systems and weighted homogeneous curves}

\thanks{The authors thank Pierre Deligne \cite{De} for his valuable remarks on \cite{Si}.}

\begin{document}

\begin{abstract} 
{We introduce a notion of rigid local system on the complement of a plane curve $Y$, which  relies on a canonical  
 Waldhausen decomposition of the Milnor  sphere associated to $Y$.  We show that 
when  $Y$   is   weigthed homogeneous this notion is deeply related to  the classical notion of rigidity   on the Riemann sphere.  
We construct large families of  rigid local  systems on the complement of  weighted homogeneous plane  curves and show that the corresponding     
  $\mathfrak D$-modules  are    generated by    `special' multivalued holomorphic functions.}
\end{abstract}

\maketitle

\section{Introduction} 
\label{sec:introduction}
{{}}
Let us consider the class of Fuchsian differential equations,        
\begin{equation}\label{fuchse}
\sum_{i=0}^na_i\frac{d^iu}{dx^i}=0,
\end{equation}
with singular points at $p_1,\ldots,p_r$ and $\infty$.
Set $p_0=\infty$.
Let $M_i\in\mbox{GL}_n(\C)$ be the monodromy of the sheaf of solutions $\mathcal S$ of  (\ref{fuchse}) along a loop $\gamma_i$ around $p_i$, $0\le i\le r$.
We call the conjugacy class of $M_i$ the \em local monodromy \em of $\mathcal S$ around $p_i$.
We will assume that we have  chosen the loops $\gamma_i$ in such a way that 
\begin{equation}\label{matrix}
M_0M_1\cdots M_r=\mathbf{I}_n.
\end{equation}
Since the singular points of (\ref{fuchse}) are regular, the degrees of the polynomials $a_i$ are bounded and (\ref{fuchse}) only
depends on a finite number $p$ of complex numbers.
Assume that we have normalized (\ref{fuchse}) replacing $u$ by 
$\prod_{i=1}^r(x-p_i)^{\lambda_i}u$, where the $\lambda_i$, $i=1,\ldots,r$,  are conveniently chosen complex numbers,
in order to maximize the number of eigenvalues of the matrices $M_i$, $1\le i\le r$, that are equal to $1$.
If $q$ denotes the number of `free' eigenvalues of the matrices $M_i$, $0\le i\le r$, then $p\ge q$.
If $p=q$ we say that (\ref{fuchse}) is \em free from accessory parameters \em \cite[Section 3.4.2]{Kohno}.

The fact that Riemann's hypergeometric differential equation is free from  accessory parameters was the key point that allowed
Riemann to compute its monodromy. Riemann's strategy can be described in the following way: 
(i) \em to classify a certain class of irreducible linear representations of a fundamental group \em \cite[Section 2.4]{GP}; 
(ii) \em to classify a class of differential equations with prescribed order and prescribed singularities  \em \cite{AHL};
(iii) \em to establish a bijection between these two classes of objects using ``local data''\em : \em 
the roots of the indicial equations of the differential equation at the singular points 
and the conjugacy classes of certain matrices associated to the representations, 
that describe the local monodromies of the differential equation. \em

Levelt \cite{Lev} established a vast generalization of Riemann's work following the same strategy.

{{}}

The theory of rigid local systems is a very ambitious reformulation, in the framework  of sheaf cohomology, of the notion of
Fuchsian differential equation free from accessory parameters  (see \cite{Kt,Simp}).

Let $\mathcal L$ be a local system on $\mathbb{P}^1\setminus \{p_0,\ldots,p_r\}$.
Let $D_i$ be  small closed disks centered at $p_i$, $0\le i\le r$.
The local system $\mathcal L$ is \em rigid \em if $\mathcal L$ is determined by its restriction to the boundaries $\partial D_i$, $0\le i\le r$.
In other words, assuming that   $\mathcal L$  is given by matrices $M_i$, $0\le i\le r$, verifying (\ref{matrix}), $\mathcal L$ is rigid if and only if 
the conjugacy classes of the matrices $M_i$ determine the simultaneous conjugacy class of the $(r+1)$-uple $(M_0,\ldots,M_r)$.
A Fuchsian differential equation is free from accessory parameters if and only if its sheaf of solutions is rigid (cf. \cite{Kt}). 
{This type of equations have been intensively  studied by several authors, notably by
Haraoka and   Yokoyama  (see, for instance,    \cite{Ha2,Ha3,HaYo,Yo}).}

Sato, Kashiwara, Kimura and Oshima \cite[Theorem 8.6 and Remark 8.7]{SKKO}
introduced a very interesting higher dimensional generalization of the notion of accessory parameters.
Let $Y=Y_{n,k,\ell}$ and  $Y^*=Y^*_{n,k,\ell}$ be  weighted homogeneous plane curves defined, respectively, by the equations 
\begin{equation}\label{w1}
\prod^\ell_{i=1}(y^k-C_ix^n)=0
\end{equation}
and
\begin{equation}\label{w2}
y\prod^\ell_{i=1}(y^k-C_ix^n)=0,
\end{equation}
where $C_i$, $i=1,\ldots,\ell$,  are pairwise distinct nonzero complex numbers and $n,k,\ell$ are positive integers such that 
$n>k$ with  $n$ and $k$  coprime integers.

\begin{theorem}\label{SKKO-th}
Set $X=\mathbb C^2$ and $\vartheta=kx\partial_x+ny\partial_y$.
Let $\mathfrak{M}$ be the germ at the origin of a simple holonomic $\mathfrak{D}_X$-module with characteristic variety    
$T^*_{Y}X\cup T^*_XX$. There are complex numbers $\lambda,\lambda_{r,s}$, $r,s\ge 0$, $ns-k(\ell+1)\le r \le (n-k)s-1$, and 
$\widehat{C}_i$, $1\le i\le \ell$, such that $\mathfrak{M}$ is isomorphic to
the system 
\begin{equation}\label{skko}
(\vartheta-\lambda)u=Pu=0,
\end{equation}
where
$P=\prod_{i=1}^\ell(\partial^k_x+\widehat{C}_ix^{n-k}\partial _y^k)+\sum_{r,s}\lambda_{r,s}x^r\partial_x^{r-ns+k(\ell+1)}\partial_y^{k s}$. 
\end{theorem}
The complex numbers $\widehat {C}_i$ are determined by the complex numbers $C_i$, $1\le i\le \ell$.
Let  $p$ be the number of parameters $\lambda,\lambda_{r,s}$ that occur in the system (\ref{skko}). 
Let $M_i$ be the local monodromy of the sheaf of solutions of $\mathfrak M$ around the regular part of $\{y^k-C_ix^n=0  \}$, $i=1,\ldots,\ell$. 
Then $M_i$ is a pseudo-reflection: 
it has  one and only one eigenvalue distinct  from $1$, with multiplicity $1$. 
Hence the number of `free' eigenvalues of the matrices $M_i$ equals $\ell$. These eigenvalues are determined by the orders of 
$\mathfrak{M}$  at the regular part of irreducible components of its characteristic variety. 
The authors of \cite{SKKO} showed that $p\ge \ell$ for each system of type (\ref{skko}). 
{Moreover, they showed that  $p= \ell$, i.e., that  $\mathfrak{M}$ is \em free from  accessory parameters\em ,     if and only if 
$\ell\leq 2$  when $\min\{k,n-k\}=1$  and $\ell=1$,  otherwise. }
Systems without accessory parameters are supposed to be interesting objects.
It is somehow  disappointing that there are so few systems with this property.
We can also find in \cite{SKKO} a similar discussion when  we replace $Y$ by $Y^*$.

{In \cite{NeSi} the authors extended  the classification  Theorem \ref{SKKO-th} to the multiplicity one case when  $Y$ is an  irreducible cusp.  
They proved  that already in this case the $\mathfrak{D}$-modules
are  no longer determined by  their   orders  at  the regular part of the  characteristic variety. 
Some (finite) additional data at the singular point has also to be considered.   
 Nevertheless, this  additional  data can 
be recovered by taking  the restriction of 
the $\mathfrak{D}$-module to certain  non-characteristic divisors passing through  the singular point of its ramification  locus.   
  We propose here  a refinement of the  notion of local data  to encompass  these local  monodromies. }

There is a close relationship between  the  systems of  PDEs of the previous two paragraphs  and   Fuchsian differential equations  on the Riemann sphere. 

Set $v(x,y)=y^{-\lambda/n}u(x,y)$. Since $\vartheta v=0$, $v$ is constant along the integral curves of $\vartheta$. 
Hence there is a multivalued holomorphic function $\varphi$ on the Riemann sphere, ramified along $0,\infty, C_i,1\le i\le \ell$, such that
$u(x,y)=y^{\lambda/n}\varphi(y^k/x^n)$. Moreover, there is a Fuchsian differential operator $H$ such that
\begin{equation}
y^{-\lambda/n}Py^{\lambda/n}\varphi (y^k/x^n)=(H\varphi)(y^k/x^n).
\end{equation}
When  $\ell=1$, the system is free from acessory parameters and
$H$ is a generalized hypergeometric differential operator  in the terminology of \cite{Lev} (see also \cite[Theorem 5.3]{NeSi}).
If in addition  $k=2$, $H$ is the   Riemann's hypergeometric differential operator.

There is a dictionary between the category of  systems  of linear differential equations on 
$\mathbb P^1$ with singular regular points at $p_0,\ldots,p_r$ and the category  of local systems on $\mathbb P^1\setminus\{p_0,\ldots,p_r\}$. 
There is a similar dictionary between a certain category of regular holonomic systems  
with solutions ramified along a singular hypersurface $Y$ 
and a certain category of  local systems on the complement of $Y$, that we call {\em Pochhammer} local systems (see Section \ref{Pochh}). 
If $L$ is a generic line transversal to $Y$, the monodromy of the restriction of a Pochhammer local
system to $L\setminus Y$  define a Pochhammer tuple (cf. \cite{DeRe}).

The main purposes of this paper   are: 
to introduce the notion of rigid local system on the complement of a  plane curve $Y$; 
to show that  when $Y$ is weighted homogeneous  our notion of rigidity is deeply related to the classical notion of rigidity on the Riemann sphere;
to show that if we replace simple holonomic system by holonomic system of multiplicity one  in the problem considered in \cite{SKKO}, 
we can prove the existence of rigid holonomic systems with solutions ramified along an weighted homogeneous curve and prescribed 
local data; to establish the foundations of a local theory of special functions on several complex variables (see Theorem \ref{speciald}, 
Definition \ref{specialf} and  Theorem \ref{mourinho}).

This paper is essentially the study of a certain class of irreducible linear representations of certain discrete groups, 
the fundamental groups of weighted homogeneous plane curves, written in the language of local systems. 
In Section \ref{planecurves} we recall some  definitions on the topology of plane curves that are necessary 
to introduce the definition of rigid local system.
 {In Section \ref{localsystems} we introduce the notion of rigid local system on the complement of a plane curve using a
canonical  Waldhausen decomposition of the Milnor sphere associated to  the curve.  When the curve is weighted homogeneous 
 we show that this notion is closely related  to  the classical definition  of rigidity on the Riemann sphere.  Roughly speaking,   it  replaces the boundaries of small disks
at the singular  points in the Riemann sphere by the boundaries of tubular neighbourhoods of the irreducible components of  the link of the curve in the Milnor sphere.}
In Section \ref{Pochh} we discuss the relation between Pochhammer local systems and Pochhammer tuples. In sections \ref{NR_turbine}, 
\ref{recrec} and \ref{WH} we construct rigid Pochhammer local systems on the complement of weighted homogeneous curves, 
 solving convenient Deligne-Simpson's type problems.  
 In Section \ref{dmodules} 
we relate the notion of Pochhammer  local system with the theory of $\mathfrak D$-modules through the Riemann-Hilbert correspondence. 
We prove in particular  that the existence of a  Pochhammer  
system on the complement of a  weighted homogeneous plane curve $Y$ implies 
the existence of a  special multivalued holomorphic function that ramifies along $Y$.

{This is the first of a series of papers on rigid local systems on the complement of  hypersurfaces. 
Our next purpose is to extend step (i) of Riemann's program referred before   to the case of an arbitrary plane curve.
We can find in \cite{Si,Si1} a construction of Pochhammer local systems on the complement of an irreducible plane curve following the same ideas. 
Most of these local systems are rigid. Steps (ii) and (iii) will be further  developped in another forthcoming paper, at least for the weighted homogeneous case.
These steps were already partially accomplished in the irreducible weighted homogeneous case (cf. \cite{NeSi}). 
}


\section{Geometric settings}\label{planecurves}

\subsection{Weighted homogeneous curves and turbines}

\vspace{1ex}

Let $n,k$ be coprime positive integers such that  $n>k$.    Let $\ell$ be a positive integer.
Following the terminology of the previous section  we set, 
\[
Y_0=\{y=0\},\qquad Y_i=\{(x,y) :y^k-C_ix^n =0\}, \qquad 1\le i\le \ell,
\]
with the $C_i$'s pairwise distinct nonzero complex numbers. Choose  $\rho>0$ small enough such that
the `fat curves', 
\begin{equation}\label{fat}
\widetilde Y_0=\{(x,y):|y|\leq \rho|x|\} ,
\quad \widetilde{Y}_i=\{(t^k,\lambda t^n): |\lambda^k-C_i|\leq\rho,\:t\in\C\}, ~1\le i\le \ell, 
\end{equation}
do not intersect each other outside a  fixed small neighborhood  of the origin.
For  $\varepsilon, \delta>0$ and  $0\le i\le \ell$ set,  
\[K_i:=Y_i\cap \partial (D_\varepsilon\times D_\delta ), \qquad N_i:=\widetilde{Y}_i\cap \partial (D_\varepsilon\times D_\delta).
\]
The topological $3$-dimensional sphere $\partial (D_\varepsilon\times D_\delta)$ is the union of  the  solid tori
$N_\infty:=\partial D_\varepsilon\times D_\delta$ and  $N^\infty:=D_\varepsilon\times \partial D_\delta$,  pieced together along their common boundary 
$\partial D_\varepsilon\times \partial D_\delta$. 
The knots $K_i$,  $1\le i\le \ell$,  are torus knots of type $(n,k)$ and  the knot $K_0$ is the  trivial knot.  For convenient choices of $\delta\gg \varepsilon>0$,  
the tubular neighbourhoods  $N_i$ of $K_i$, 
 $0\le i\le \ell$,   are pairwise  disjoint and  contained in the interior of $N_\infty$.  Set 
 \[ {{T}}_{n,k,\ell}=N_\infty\setminus\cup_{i=1}^\ell \mbox{int}(N_i), \qquad {{T}}^*_{n,k,\ell}=N_\infty\setminus\cup_{i=0}^\ell \mbox{int}(N_i). \]
 We call a topological space homeomorphic to ${{T}}_{n,k,\ell}$ [${{T}}^*_{n,k,\ell}$] a \em  turbine \em without  [with]  shaft and parameters $n,k,\ell$.
The turbines   ${{T}}_{n,k,\ell}$, ${{T}}^*_{n,k,\ell}$   are, respectively,  retracts  by deformation of the   complement  of the  plane curves,
\begin{equation}\label{eq1}
x{\prod_{i=1}^\ell}(y^k-C_ix^n)=0,
\end{equation}
\begin{equation}\label{eq2}
xy{\prod_{i=1}^\ell}(y^k-C_ix^n)=0.
\end{equation}
 
Given an weighted homogeneous plane curve $Y$ we   shall denote by $T_Y$ the  turbine associated to $Y$  by the construction  above. 

\vspace{1ex}

For  each $i=0,\ldots,\ell,\infty$, let $\alpha_i,\beta_i$ be the homotopy classes of positively oriented simple   closed curves  on $\partial N_i$ such that $\alpha_i\sim 0$  
and  $\beta_i\sim K_i$ in $H_1(N_i)$,   $\ell(K_i,\alpha_i)=1$ and   $\ell(K_i,\beta_i)=0$,  where
$\ell(\cdot,\cdot)$ denotes the linking number inside the oriented homology  3-sphere $\partial(D_\varepsilon\times D_\delta)$. 
For each  $i=0,\ldots,\ell,\infty$, the pair   $\alpha_i,\beta_i$  generates  $\pi_1(\partial N_i)$ and  is 
unique up to isotopy.   We call  $\alpha_i$  [$\beta_i$],    the {\em standard meridian}  [{\em   parallel}] of the torus $\partial N_i$.

For   $i=0,\infty$, set $\delta_i=\alpha_i^n\beta_i^k$ and $\omega_i=\alpha_i^s\beta_i^r$, 
where  $r,s$ are  integers such that $rn=ks+1$.  Then $\alpha_i=\delta_i^r\omega_i^{-k}$ and  $\beta_i=\delta_i^{-s}\omega_i^{n}$, $i=0,\infty$.
By the results of  \cite[Lemma 2.3]{NeSi1}, we have $\delta_0=\delta_\infty=\alpha_i^{nk}\beta_i$, $i=1,\ldots,\ell$, in  $\pi_1(T^*_{n,k,\ell})$,   which  we denoted  by  $\delta$.
Thus,  $\omega_i,\delta$  generate $\pi_1(\partial N_i)$ $i=0,\infty$ and   $\alpha_i,\delta$  generate
$\pi_1(\partial N_i)$ $i=1,\ldots,\ell$. Moreover,  next result holds.

\begin{theorem}\label{gerad_proj} The following hold.
\begin{enumerate}
\item[$(a)$] $\pi_1(T^*_{n,k,\ell})=\langle\alpha_1,\ldots,\alpha_\ell,\omega_0,\omega_\infty
:\alpha_1\cdots\alpha_\ell\omega_\infty=\omega_0\rangle\times\langle \delta\rangle$.
\item[$(b)$]  $\pi_1(T_{n,k,\ell})$ is  isomorphic to the quotient group  of  $\pi_1({{T}}^*_{n,k,\ell})$ by the relation $\omega_0^k=\delta^r$.
\item[$(c)$] The  fundamental group of the complement of  the curve $Y^*_{n,k,\ell}$  {\em[}$Y_{n,k,\ell}${\em]} 
is isomorphic to the quotient group   of $\pi_1({{T}}^*_{n,k,\ell})$ {\em[}$\pi_1({{T}}_{n,k,\ell})${\em]} by the relation  $\omega_\infty^n=\delta^s$.
\end{enumerate}
\end{theorem}
{\bf Proof.}  $(a)$ follows from  \cite[Lemma 2.3]{NeSi1}. The remaining assertions follow from $(a)$  and obvious geometric considerations.\quad$\Box$

\begin{remark}\label{gerar_proj_rmk}
The presentations given in  Theorem \ref{gerad_proj} are independent  of the choice of the pair of integers $r,s$ verifying $rn=sk+1$. Actually, 
if $r',s'\in\mathbb{Z}$ is another pair  verifying $r'n=s'k+1$, there is $t\in\mathbb{Z}$ such that $r'=r+tk$ and $s'=s+tn$, and 
it is enough to replace $\omega_i$,   $i=0,\infty$,  by $\omega'_i=\omega_i\delta^t$.
\end{remark}

From Theorem \ref{gerad_proj}   we derive easily  the following well-known  presentation of the  local fundamental group of a cusp.
\begin{corollary}\label{gerad_proj-curve} 
The fundamental group of the  complement of the cusp $y^k=x^n$,   is given by
\[ \pi_1(Y_{n,k,1})=\langle \alpha_\infty,\beta_0 \::\:\alpha_\infty^n=\beta_0^k\rangle.   \]
\end{corollary}
We conclude with another corollary of Theorem \ref{gerad_proj} that will be very useful  in  Section \ref{recrec}. 
\begin{corollary}\label{gerad_proj_corol}
We have the    relation in $\pi_1(T^*_{n,k,\ell})$,  
\[Ê\alpha_\infty=g_0(g_{0,1}\,\cdots\,g_{0,\ell})\,( g_{1,1}\,\cdots\, g_{1,\ell})\,\cdots\,(g_{k-1,1}\,\cdots\, g_{k-1,\ell}),  \]
where $g_0=\alpha_0$, $g_{j,i}=\omega_0^{k-1-j}\alpha_i\omega_0^{1-k+j}$, $j=0,\ldots,k-1$ and $i=1,\ldots,\ell$. 
Moreover by letting $g_0=1$ in relation above, we obtain the  analogous  expression for $\alpha_\infty$  in $\pi_1(T_{n,k,\ell})$. 
\end{corollary}
\noindent{\bf Proof.} By Theorem \ref{gerad_proj} along with relations $\alpha_i=\delta^r\omega_i^{-k}$, $i=0,\infty$ we get,
\begin{eqnarray*}
\alpha_\infty    &=&\delta^r\omega_\infty^{-k}\\ 
	&=&\delta^r(\omega_0^{-1}\alpha_1\cdots\alpha_\ell)^k\\
           &=&\alpha_0\omega_0^k(\omega_0^{-1}\alpha_1\cdots\alpha_\ell)^k\\
           &=&\alpha_0\omega_0^{k-1}(\alpha_1\cdots\alpha_\ell)(\omega_0^{1-k}\omega_0^{k-2})(\alpha_1\cdots\alpha_\ell)(\omega_0^{2-k}
               \omega_0^{k-3})\cdots\omega_0^{-1}(\alpha_1\cdots\alpha_\ell)\\
           &=&\alpha_0(\omega_0^{k-1}\alpha_1\cdots\alpha_\ell\omega_0^{1-k})({\omega}^{k-2}\alpha_1\cdots\alpha_\ell{\omega}^{2-k})\cdots
              (\alpha_1\cdots\alpha_\ell).\quad\Box
\end{eqnarray*}

\subsection{Waldhausen decomposition of the Milnor sphere associated to a  plane curve}

\vspace{1ex}

Let $V$ be  a 3-dimensional manifold. We say that   $V$ is a {\em Seifert manifold}  if $V$  is  decomposed   into circles, the {\em fibers},
such that  each fiber   has a  tubular neighborhood difeomorphic  (preserving fibers) to a {\em standard fibered  torus}.

Let $M^3$ be an oriented connected compact 3-dimensional  manifold. 
A {\em Waldhausen decomposition} of $M^3$   is a finite partition  $M^3=\bigsqcup_k \mathcal{T}_k\bigsqcup_jV_j$ where  $(\mathcal{T}_k)$ is  family of  
2-tori and   $(V_j)$ a family  of     Seifert manifolds.  
We say that the Waldhausen decomposition    is  {\em minimal}  if the number of 2-tori $(\mathcal{T}_k)$ is minimum. 
We say that  the Waldhausen decomposition  is {\em adapted} to a given  1-dimensional submanifold of $M^3$ if each irreducible component
of this  manifold is a   fiber of  one of the 
Seifert manifolds $V_j$ (see  \cite{Le1}).

Let $(Y,0)$ be a germ of a plane curve defined by a reduced analytic function $f(x,y)\in\C\{x,y\}$. Assume moreover that the tangent cone of $Y$ is 
transversal to $\{x=0\}$. 
For $0<\varepsilon\ll 1$, the 3-dimensional sphere $S^3_\varepsilon=\{ (x,y): |x|^2+|y|^2=\varepsilon^2\}$  
intersects $Y$  transversally into an (algebraic) link $\mathbf{L}$.
By Milnor's cone Theorem,  the topology of the complement of $Y$  is determined by the topology of the 3-dimensional compact connected  manifold 
$S^3_\varepsilon\setminus\mbox{int}(N(\mathbf{L}))$, where $N(\mathbf{L})$ is a sufficiently  thin closed tubular neighbourhood of $\mathbf{L}$. 

We decompose the  {\em Milnor sphere} $S^3_\varepsilon$ into   
 the union of two solid tori $N_\infty$ and $N^\infty$ pieced together  along their common boundary  $\partial N_\infty=\partial N^\infty$ 
such that  $N(\mathbf{L})\subset \mbox{int}(N_\infty)$.  Let    $K^\infty$ be the core of $N^\infty$.

\begin{theorem}\label{minimal}\em \cite[Theorem 4.5.1]{Le1} \em
There is a minimal Waldhausen decomposition  of $S^3_\varepsilon$ adapted to the components of the 
 link $\mathbf{L}\cup K^\infty$. 
 Moreover,  this decomposition  is  unique {\em(}up to isotopy{\em)}.
\end{theorem}

\begin{remark}\label{pararigidezlocal}
We can find in  \cite{NeSi1} a simple construction of a    minimal  Waldhausen decomposition of $S^3_\varepsilon$ 
determined by a  family of tori $(\mathcal{T}_k)$ that  includes the torus $\partial N_\infty$,  
such that the closure of each connected component of $N_\infty\setminus\bigcup_k\mathcal{T}_k$ is a turbine.  
Furthermore, this decomposition is described by a   tree  similar to the Eggers tree of $Y$.   
\end{remark}



\section{Rigid local systems}\label{localsystems}

\subsection{Review on rigid local systems on the Riemann sphere}

\vspace{1ex}

Set  $U=\mathbb{P}^1\setminus S$ where $S=\{0,C_1,\ldots,C_\ell,\infty\}$ is a  subset of $\mathbb{P}^1$ with $\ell+2$ elements,  and let
 $j:U\hookrightarrow \mathbb{P}^1$ be the inclusion map. Set $C_0=0$ and $C_\infty=\infty$. 
Let  $\widetilde{D}_i\subset \mathbb{P}^1$, $i=0,1,\ldots,\ell,\infty$, be  disjoint closed disks  of nonzero  radii  centered at $C_i$, $i=0,1,\ldots,\ell,\infty$, respectively.

Let $\mathcal{F}$  be  local system on $U$.   We call {\em local data}  of $\mathcal{F}$  to the collection of 
local systems  $\mathcal{F}|_{\partial\widetilde{D}_i}$, $i=0,1,\ldots,\ell,\infty$. We say that $\mathcal{F}$  is  {\em   rigid}  if  $\mathcal{F}$ 
is determined by  its local data (up to isomorphism).
We say that $\mathcal F$ is {\em semi-rigid} if there are only a finite number of isomorphism classes of  local systems on $U$
with local data  isomorphic to the local data of $\mathcal{F}$.

The   {\em  index of rigidity}  of   $\mathcal{F}$ in $U$  is defined as 
$\mbox{rig}(\mathcal{F},U)=\chi(\mathbb{P}^1,j_*\mbox{End}(\mathcal{F}))$. 
Using  the  Euler-Poincar\'e formula we derive the formula,
\[ \mbox{rig}(\mathcal{F},U)=
\sum_{i\in\{0,1,\ldots,\ell,\infty\}}\dim \Gamma(\partial\widetilde{D}_i,\mbox{End}(\mathcal{F}))-\ell\:(\mbox{rank}\:\mathcal{F})^2. 
\]

\begin{theorem}\em\cite[Theorem 1.1.2]{Kt}\em \label{Katz_result}
An irreducible local system $\mathcal{F}$ on $U$ is rigid  if and only if is semi-rigid if and only if  $\mbox{rig}(\mathcal{F},U)=2$.
\end{theorem}

\subsection{Rigid local systems on the complement of a plane curve}

\vspace{1ex}

\begin{definition}\label{rigid-def}
Let $\mathcal{L}$ be a  local system defined on the complement of the  germ (at the origin) of a  plane curve $Y$ with tangent cone transversal to $\{x=0\}$. 
Let $S^3_\varepsilon=\bigsqcup_k\mathcal{T}_k\bigsqcup_j V_j$ be a  minimal Waldhausen decomposition of a Milnor sphere 
$S^3_\varepsilon$ in the conditions of Theorem \ref{minimal}.
We call   {\em local data} of $\mathcal{L}$ to the collection of local systems $({\mathcal L}|_{ \mathcal{T}_k})$. 
We say that $\mathcal L$ is {\em  $\mathcal{W}$-rigid}  if the isomorphism class of $\mathcal L$ is determined by the local data of $\mathcal{L}$.
We say that $\mathcal L$ is {\em  $\mathcal{W}$-semi-rigid}    if there are only a finite number of isomorphism classes of  local systems  on the complement of  $Y$
with local data   isomorphic to the local data of $\mathcal{L}$.
\end{definition}

By Theorem \ref{minimal} the minimal Waldhausen decomposition   is unique modulo an isotopy. Hence
the definition of local data does not depend on  the choice of   this   decomposition.

We shall refer  to a $\mathcal{W}$-rigid [$\mathcal{W}$-semi-rigid] local system  simply as  rigid [semi-rigid].

\begin{definition}\label{rigid-def-turbine} Let  $\mathcal{L}$ be a local system on a turbine ${T}$.  
We call {\em local data} of $\mathcal{L}$ to $\mathcal{L}|_{\partial T}$. 
We say that $\mathcal L$ is {\em  rigid}   if  the isomorphism class of $\mathcal L$ is determined by the local data $\mathcal{L}$.
We say that $\mathcal L$ is {\em semi-rigid} if there are only a finite number of isomorphism classes of  local systems on $T$
with local data  isomorphic to the local data of $\mathcal{L}$.
\end{definition}

The  connected components of  the boundary of a turbine $T$ are  tori, say  $\mathcal{T}_1,\ldots,\mathcal{T}_t$, and henceforth  
  the    local data of  $\mathcal{L}$ is the  collection of    simultaneous conjugacy classes  of the  pairs $(\varrho(a_i),\varrho(b_i))$, $i=1,\ldots,t$,
with   $a_i,b_i$   generators of  $H_1(\mathcal{T}_i)$ and $\varrho:\pi_1(T)\to\mbox{GL}(m)$  the monodromy representation of $\mathcal{L}$.


In view of  Remark \ref{pararigidezlocal}  we have immediately the following result.
\begin{lemma}\label{restr-turbine}
Let $Y$ be an weighted homogeneous   plane   curve with tangent cone $\{y=0\}$ and let  $T_Y$ be  a turbine associated to $Y$. 
Let $\mathcal{L}$ be a local system defined on the complement of $Y$.  Then   $\mathcal{L}$ is rigid if and only if   $\mathcal{L}|_{T_Y}$  is rigid.  
A similar conclusion holds if we replace rigid by semi-rigid.
\end{lemma}

Let $Y^\circ_{n,k,\ell}$ be the curve defined by (\ref{eq2}).  
Given $\varepsilon\in\C^*$, let  $\widehat{\varepsilon}$  be the one-dimensional representation  of
$\pi_1(\mathbb C^2\setminus Y^\circ_{n,k,\ell})\simeq\pi_1(T^*_{n,k,\ell})$ given by  $\widehat{\varepsilon}(\delta)=\varepsilon$ and 
$\widehat{\varepsilon}(\alpha_i)=\widehat{\varepsilon}(\omega_0)=1$, $1\le i\le \ell$.
Let  $\mathcal K_\varepsilon$  be the rank one local system  on $\mathbb C^2\setminus Y^\circ_{n,k,\ell}$ whose monodromy representation is $\widehat{\varepsilon}$. 

Set $S=\{C_0,C_1,\ldots,C_\ell, C_\infty\}$ and consider the map
$\gamma:\mathbb C^2\setminus Y^\circ_{n,k,\ell}\to\mathbb{P}^1\setminus S$ given  by $\gamma(x,y)=(y^k:x^n)$.   The homotopy class of  each  fiber
of $\gamma$   is given by  $\delta$.
Let $\Sigma^\varepsilon$ be the category of local systems on $\mathbb C^2\setminus Y^\circ_{n,k,\ell}$ with scalar monodromy $\varepsilon\cdot\mathbf{Id}$ along
the fibers of $\gamma$ and let $\mathcal L\mapsto \underline{\mathcal L}$ be the functor that associates to an object of $\Sigma^\varepsilon$ the unique 
 local system $\underline{\mathcal L}$ on $\mathbb P^1\setminus S$ such that 
 $\mathcal K^{\otimes -1}_\varepsilon\otimes \mathcal L\simeq \gamma^{-1}\underline{\mathcal L}$.

\begin{lemma}\label{equivalen}
The functor 
 $\mathcal{L}\mapsto \underline{\mathcal{L}}$ 
is an equivalence of categories such that  $\mathcal{L}$ is irreducible {\em[}rigid, semi-rigid{\em]}  
 if and only if $\underline{\mathcal L}$ is irreducible  {\em[}rigid, semi-rigid{\em]}.
\end{lemma}
{\bf Proof.} Let $\mathcal L$ be an object  of $\Sigma^\varepsilon$. 
The monodromy of $\mathcal K_\varepsilon^{\otimes -1}\otimes \mathcal L$ is trivial along the fibers of $\gamma$. Hence
$\underline{\mathcal L}:=\gamma_*(\mathcal K_\varepsilon^{\otimes -1}\otimes\mathcal{L})$  is a local system   on $\mathbb P^1\setminus S$.
Conversely, if $\underline{\mathcal L}$ is a local system on $\mathbb{P}^1\setminus S$, $\gamma^{-1}\underline{\mathcal{L}}$ is a local system on 
the complement of $Y^\circ_{n,k,\ell}$ with   trivial  monodromy   along  the fibers of $\gamma$. Thus
 $\mathcal K_\varepsilon\otimes \gamma^{-1} \underline{\mathcal L}$
 has scalar  monodromy  $\varepsilon$ along the fibers of $\gamma$. Moreover,  
 \[ \mathcal K_\varepsilon\otimes \gamma^{-1} \left(\gamma_*(\mathcal K_\varepsilon^{\otimes -1}\otimes\mathcal{L})\right)\simeq \mathcal{L},\qquad
 \gamma_*\left(\mathcal K_\varepsilon^{\otimes -1}
 \otimes\left(\mathcal K_\varepsilon\otimes \gamma^{-1} \underline{\mathcal L}\right)\right)\simeq \underline{\mathcal L}.\] 


By Theorem \ref{gerad_proj},   $\pi_1(\mathbb C^2\setminus Y^\circ_{n,k,\ell})\simeq\pi_1(\mathbb P^1\setminus S)\times\langle\delta\rangle$,
where $\langle\delta\rangle$ is the infinite cyclic group generated by $\delta$. 
Hence   $\mathcal{L}$ is irreducible  if and only if $\underline{\mathcal L}$ is irreducible.

Let   $\bar\gamma$ be the restriction of  $\gamma$ to  $\partial (D_\varepsilon\times D_\delta)\setminus  Y^\circ_{n,k,\ell}\equiv T^*_{n,k,\ell}$. 
There are disjoint closed disks $\widetilde{D}_i\subset\mathbb P^1$ 
with centers $C_i$, $i=0,1,\ldots,\ell,\infty$, such that
 $\overline{\gamma}^{-1}(\partial \widetilde D_i)=\partial N_i$, $i=0,\ldots,\ell$, and 
$\overline{\gamma}^{-1}(\partial \widetilde D_\infty)=\partial N_\infty=\partial N^\infty$. Thus  a local system $\mathcal L\in\Sigma^\varepsilon$ 
is determined by  $\mathcal L|_{\partial {{T}}^*_{n,k,\ell}}$ if and only if 
$(\mathcal K_\varepsilon^{\otimes -1}\otimes \mathcal L)\in \Sigma^1$ is  determined by 
$(\mathcal K_\varepsilon^{\otimes -1}\otimes \mathcal L)|_{\partial {{T}}^*_{n,k,\ell}}$   if and only if $\underline{\mathcal L}$ is determined by its 
restriction to  $\bigcup_{i=0,\ldots,\ell,\infty}\partial \widetilde{D}_i$.   Hence $\mathcal{L}$ is rigid [semi-rigid] 
 if and only if  $\underline{\mathcal L}$ is rigid  [semi-rigid].
\quad$\Box$
 \vspace{2ex}

In the sequel  we   identify a  local system on $\mathbb C^2\setminus Y^\circ_{n,k,\ell}$ with its restriction  to  $T^* ={{T}}^*_{n,k,\ell}$.


\begin{theorem}\label{index-rigidity}
If    $\mathcal{L}^*$ is an irreducible local system on   $T^*$,   $\mathcal{L}^*\in\Sigma^\varepsilon$ for some $\varepsilon\neq 0$ and 
\[ \mbox{rig}(\underline{\mathcal{L}}^*,\mathbb{P}^1\setminus S)=\sum_{i\in\{0,1,\ldots,\ell,\infty\}}\dim\Gamma(\partial N_i,\mbox{End}(\mathcal{L}^*))  -\ell(\hbox{\rm rank}\:\mathcal L^*)^2. \]
Moreover, $\mathcal{L}^*$ is rigid if and only if $\mbox{rig}(\underline{\mathcal{L}}^*,\mathbb{P}^1\setminus S)=2$. 
 \end{theorem}
 {\bf Proof.}
 Keep the notations of the proof of Lemma \ref{equivalen}. 
Since $\delta $ is in the center of  $\pi_1(T^*)$ and $\mathcal{L}^*$  is irreducible, 
$\mathcal{L}^*$ has    scalar monodromy along the fibers of $\overline{\gamma}$.  Hence $\mathcal{L}^*\in\Sigma^\varepsilon$ for some $\varepsilon\neq 0$.
By Lemma \ref{equivalen}  there is  an irreducible local system  $\underline{{\mathcal{L}}}$ 
 on $\mathbb{P}^1\setminus S$  such that $\mathcal{L}^*\simeq(\mathcal{K}_\varepsilon)|_{T^*}\otimes \overline{\gamma}^{-1}\underline{\mathcal L}$.
 Since $\mbox{rank}(\mathcal{K}_\varepsilon)=1$, $\mbox{rank}(\mathcal{L}^*) =\mbox{rank}(\underline{\mathcal L})$.  Since 
  $\overline{\gamma}^{-1}(\partial \widetilde{D}_i)=\partial N_i$  for  $i=0,1,\ldots,\ell,\infty$,   $\dim\Gamma(\partial N_i, \mbox{End}(\mathcal{L}^*))=\dim\Gamma(\partial \widetilde{D}_i,\mbox{End}(\underline{\mathcal{L}}))$. 
The proof now  follows from     Theorem \ref{Katz_result} and   Lemma \ref{equivalen}.\quad$\Box$
 \vspace{2ex}

We  denote $\mbox{rig}(\mathcal{L}^*,T^*)=\mbox{rig}(\underline{\mathcal{L}}^*,\mathbb{P}^1\setminus S)$, which we refer as 
 the  {\em index of rigidity} of  $\mathcal{L}^*$.


For the turbine without shaft  we have a weaker result. 

\begin{theorem}\label{rig_without}
Let $\mathcal{L}$ be a local system on  $T$. The following holds:
\begin{enumerate}
\item[$(i)$] $\mathcal{L}$ is irreducible if and only if  $\mathcal{L}|_{T^*}$ is  irreducible. 
\item[$(ii)$]\label{rigid-without} If $\mathcal{L}$ is  rigid,   $\mathcal{L}|_{T^*}$ is rigid.
\item[$(iii)$]\label{rigid-without-index} If  $\mathcal{L}$ is irreducible and rigid,  $\mbox{rig}(\mathcal{L}|_{T^*},T^*)=2$. 
\end{enumerate}

\end{theorem}

{\bf Proof.}~
$(i)$ follows from the fact that  $\pi_1(T)$ is a quotient group of  $\pi_1(T^*)$ (cf. Theorem \ref{gerad_proj}).
Assume that $\mathcal{L}$ is rigid.  Let $\mathcal{F}^*$ be a local system defined in  $T^*$ with the same local data
of ${\mathcal{L}}|_{T^*}$, that is,     $\mathcal{L}|_{\partial N_i}\simeq\mathcal{F}^*|_{\partial N_i}$, $i=0,\ldots,\ell,\infty$. Since 
$\mathcal{L}|_{\partial N_0}\simeq\mathcal{F}^*|_{\partial N_0}$, we can glue  $\mathcal{F}^*$ and
$\mathcal{L}|_{N_0}$ along $\partial N_0$ and  construct a local system $\mathcal{F}$ 
on  $T$  with the same local data of $\mathcal{L}$. 
Since  $\mathcal{L}$ is rigid,  $\mathcal{L}\simeq\mathcal{F}$.  Hence $\mathcal{L}|_{T^*}\simeq\mathcal{F}^*$.

$(iii)$ follows from  Theorem \ref{index-rigidity} taking into account $(i)$ and $(ii)$.\quad$\Box$
 \vspace{2ex}



\section{Pochhammer local systems}\label{Pochh}
 
 \subsection{Pochhammer local systems and Pochhammer tuples}
 
 \vspace{1ex}
 

Let $A=\{a_1,\ldots,a_n\}$ be a finite subset of a simply connected open subset $\Omega$ of $\mathbb C$. 
Set $U=\Omega\setminus A$.  Let $j:U\hookrightarrow\Omega$ be the open inclusion. 
Let $\mathcal L$ be a local system on $U$ and  set $V=\mathcal L_a$,  for some fixed point  $a\in U$. 
Consider loops  $\gamma_i\in\pi_1(U,a)$, $1\le i\le n$, verifying
\[
\frac{1}{2\pi\sqrt{-1}}\oint\limits_{\gamma_i}\frac{dx}{x-a_j}=\delta_{i,j},\qquad  j=1,\:\ldots,\:n.
\]
Let $\varrho:\pi_1(U,a)\to\mbox{GL}(V)$ be the monodromy  representation of  $\mathcal L$.
For each $\gamma\in \pi_1(U,a)$ set $V^\gamma=\ker(\varrho(\gamma)-\mathbf{1}_V)$.

\begin{lemma}\label{lemmag}{\em \cite{N,Si}\em}
Let $\mathcal L$ be a local system on $U$.
The following statements are equivalent:
\begin{enumerate}
\item[$(a)$]  There is a unique decomposition $V=\oplus_{i=1}^n W_i$ such that $V^{\gamma_i}=\oplus_{j\not=i}W_j$;
\item[$(b)$] $\dim V=\sum_{i=1}^n\mbox{codim}\:V^{\gamma_i}$  and $\bigcap_{i=1}^nV^{\gamma_i}=\{0\}$;
\item[$(c)$] $H^*(\Omega,j_*\mathcal L)=0$.
\end{enumerate}
\end{lemma}
Note that Pochhammer local systems were called {\em hypergeometric} in  \cite{N,Si}.


\begin{definition}\label{P1-Pochhammero}
We say that   a local system  $\mathcal L$ on $U$  is   {\em Pochhammer} if it verifies 
 the equivalent conditions of   Lemma \ref{lemmag}.  
\end{definition}

 A matrix $M\in\mbox{GL}_n(\mathbb{C})$ is  a   {\em pseudo-reflection} if $\mbox{codim}\:\ker (M-\mathbf{I}_{n})=1$. 
 The determinant of a pseudo-reflection $M$   is an eigenvalue of $M$, called the   {\em special eigenvalue} of $M$. 
Unless otherwise stated, we always  assume the special eigenvalue of a pseudo-reflection  distinct from one.

Assume that    $\mathcal L$ is   Pochhammer. We call  {\em multiplicity} of $\mathcal L$ at a point $a_i\in A$  to the codimension of $V^{\gamma_i}$ and 
we denote it by  $\mbox{mult}_{a_i}(\mathcal{L})$.  If $\mbox{mult}_{a_i}(\mathcal{L})=1$, we  call  
{\em special eigenvalue} of $\mathcal{L}$ at $a_i$,  $i=1,\ldots,n$,    to
 the special eigenvalue  of the   pseudo-reflection $\varrho(\gamma_i)$.   
 We say  that   $\mathcal L$ has {\em multiplicity one}   if   $\mbox{mult}_{a_i}(\mathcal{L})=1$ for $i=1,\ldots,n$.

Assume that $U=\Pp^1\setminus\{a_1,\ldots,a_n,\infty\}$ and choose  simple loops 
$\gamma_i$ around $a_i$, $i=1,\ldots,n$, in such a way that      $\gamma_\infty:=\gamma_1\cdots\gamma_n$ 
is a simple  loop  around  $\infty$.   Denote by  $A_i$, $i=1,\ldots,n,\infty$,     the matrices  $\varrho(\gamma_i)$, $i=1,\ldots,n,\infty$,  w.r.t. some basis   of $V$. 

Next result follows from \cite[Theorem 1.1]{Ha} and Theorem \ref{Katz_result}.

\begin{theorem} \label{DSP2} 
Consider nonzero complex numbers  $\lambda_1,\ldots,\lambda_n,\eta_1,\eta_2$,   such that $\lambda_i\neq 1,\eta_1$ for all $i$, $\eta_1\neq\eta_2$,  
and $\lambda_1\cdots\lambda_n=\eta_1^{n-1}\eta_2$. For  each $i=1,\ldots,n$ set 
\begin{equation}\label{A_i-Pochhammer}
           A_i=\left(\begin{array}{cccccccc} 1 &  &&& (\lambda_1-\eta_1)\eta_1^{-1}\\ 
                                                      &\ddots&&&\vdots\\ 
                                                      && &1&(\lambda_{i-1}-\eta_1)\eta_1^{-1}\\ \\
                                                      &&&&\lambda_i\\ \\
                                                      &&&&\lambda_{i+1}-\eta_1&1\\ 
                                                      &&&&\vdots&&\ddots\\ 
                                                      &&&&\lambda_n-\eta_1&&&1
              \end{array}\right).
 \end{equation}
 Set $A_\infty:=A_1\cdots A_n$. 
Then    $(A_1,\cdots, A_n,A_\infty)$ determines  the monodromy representation of an  irreducible rigid Pochhammer  local system of multiplicity one on $U$  and  
special eigenvalue $\lambda_i$ at $a_i$, $i=1,\ldots,n$, such that $A_\infty$ is conjugated to $\eta_1\mathbf{I}_{n-1}\oplus\eta_2$.
\end{theorem}

{We call  the tuple $(\lambda_1,\ldots,\lambda_n;\eta_1,\eta_2)$ of Theorem \ref{DSP2}   the {\em numerical local data} of $\mathcal{L}$ and  
the  tuple of matrices $(A_1,\ldots,A_n,A_\infty)$  a      {\em Pochhammer tuple} (cf. \cite{DeRe}), which   determines the monodromy 
of  a Pochhammer local system of  differential equations with regular singular points  at $a_1,\ldots,a_n,\infty$ and no logarithmic  solution.
We refer to \cite{Ha} for a definition of Pochhammer system of differential equations.}

{Next result shows  that the  local systems of theorem  Theorem \ref{DSP2} exhausted  the class of rigid  irreducible 
Pochhammer  local systems of multiplicity one   on the punctured  Riemann sphere (up to isomorphism).}
 
\begin{theorem}\label{cond-gauss}
Let $\mathcal{L}$  be  an irreducible  rigid Pochhammer local system of multiplicity one  on $U$ with semi-simple monodromy around infinity. 
For each  $i=1,\ldots,n$,  let $A_i$  be the  pseudo-reflection determined by the local monodromy of $\mathcal{L}$ around  $a_i$ and let 
 $\lambda_i\neq 1$ be  its  special eigenvalue.  
There are distinct nonzero  complex numbers $\eta_1, \eta_2$ verifying  $\lambda_1,\ldots,\lambda_n\neq\eta_1$ and 
$\lambda_1\cdots\lambda_n=\eta_1^{n-1}\eta_2$,
 such that  the monodromy around infinity , $A_\infty=A_1\cdots A_n$,   is conjugated to $\eta_1\mathbf{I}_{n-1}\oplus\eta_2$.
In particular, $(A_1,\ldots,A_n,A_\infty)$ is a Pochhammer tuple. 
\end{theorem}
{\bf Proof.}
Since $\mathcal{L}$ has multiplicity one and  special eigenvalue  $\lambda_i\neq 1$ at $a_i$ all   centralizers 
$Z(A_i)=\{B\in\mbox{GL}(n):A_iB=BA_i\}$, $i=1,\ldots,n$,  have  dimension $(n-1)^2+1$.  
Since  $\mathcal{L}$ is rigid  we have by Theorem \ref{Katz_result}, 
\[  
\mbox{rig}(\mathcal{L},U)=-(n-1)n^2+n((n-1)^2+1)+\dim Z(A_\infty)=2.
\]
Hence $\dim Z(A_\infty)=(n-1)^2+1$. Let $d_1\leq d_2\leq \cdots\leq d_s$ be the multiplicities of the distinct eigenvalues  of $A_\infty$.  Clearly $s\geq 2$. 
Since  $A_\infty$ is semi-simple,  $\sum_{i=1}^sd_i=n$ and
$\sum_{i=1}^sd^2_i=\dim Z(A_\infty)=n^2+2(1-n)$. Since
\[\left(\sum_{i=1}^sd_i\right)^2= \sum_{i=1}^sd^2_i+2\sum_{i<j}d_id_j=n^2,  \] 
$\sum_{i<j}d_id_j=n-1$. In particular, $d_1(n-d_1)\leq n-1$. Since  $d_1(n-d_1)> n-1$, if  $d_1>1$, $s=2$ and $d_1=1$.  Therefore 
$A_\infty$ is  conjugated to $\eta_1\mathbf{I}_{n-1}\oplus\eta_2$
for some (nonzero) distinct  complex numbers $\eta _1\neq \eta_2$. 
The theorem now follows from the irreducibility conditions of \cite[Proposition 1.3]{Ha} and  relation $\det(A_\infty)=\det(A_1\cdots A_n)$.
\quad$\Box$\vspace{2ex}

Let $T_Y=N\setminus\bigcup_{i\in I}N_i$ be a turbine  associated to a weighted homogeneous plane curve $Y$ with
irreducible componentes $Y_i$, $i\in I$.  Set   $U=T_Y\cap L$,  with $L$ a generic line transversal to $Y$.  
Then $U$  is  a retract by deformation of     $L\setminus A$, where $A=L\cap Y$. 
\begin{definition}\label{Pochhammer-turbine}
We  say that a local system  $\mathcal F$ on $T_Y$ is    \em Pochhammer \em  if $\mathcal{F}|_{U}$  is a Pochhammer local system on  $U$.
We say  that $\mathcal F$ has {\em multiplicity} $\mu_i$ along $\partial N_i$ if   $\mathcal{F}|_{U}$  has multiplicity $\mu_i$ at some point of $Y_i\cap L$.
We say that $\mathcal{L}$ has {\em multiplicity one} on $T_Y$ if   $\mathcal{F}|_{U}$ has multiplicity one on $U$. 
\end{definition}

Next result shows that an analogue of  Theorem \ref{index-rigidity} also holds for turbines without shaft,  at least for Pochhammer local systems.

\begin{lemma}\label{matrix-infty}
Assume that  $\mathcal{L}$   is irreducible and Pochhammer on a turbine with shaft  $T_{n,k,\ell}$. 
 Then  $\mathcal{L}$ is rigid if and only if  $\mbox{rig}( \mathcal{L}|_{T^*_{n,k,\ell}}, \mathcal{L})=2$. 
\end{lemma}
\noindent{\bf Proof.}  Let  $\mathcal{K}$ be a local system on $T=T_{n,k,\ell}$  with the same local data of $\mathcal{L}$. Set
$T^*=T^*_{n,k,\ell}$  and set $\mathcal{K}^*=\mathcal{K}|_{T^*}$.
  If $\mathcal{K}|_{\partial N_0}\simeq \mathcal{L}|_{\partial N_0}$,  $\mathcal{K}^*$ and $\mathcal{L}^*$ would 
have   have the same local data on the shaft, which imply,  by the rigidity of   $\mathcal{L}^*$ that
$\mathcal{L}^*\simeq\mathcal{K}^*$. Moreover, since   $\mathcal{L}|_{\partial N_0}$ has  trivial monodromy,  
  $\mathcal{L}\simeq\mathcal{K}$. Hence it is enough to prove that   $\varrho_{\mathcal{L}}(\omega_0)$ and  $\varrho_{\mathcal{K}}(\omega_0)$ are conjugated,
  where  $\varrho_{\mathcal{H}}:\pi_1(T)\to\mbox{GL}_m(\mathbb{C})$ denote the monodromy representation of a local system $\mathcal{H}$ of rank $m$ on $T$. 

By $(i)$  $\varrho_{\mathcal{L}}(\alpha_\infty)$  does not have the eigenvalue one. Since this condition is invariant by conjugation, 
the same holds for  $\varrho_{\mathcal{K}}(\alpha_\infty)$.
By   condition $(b)$  of Lemma \ref{lemmag} $\mathcal{K}$ is a  Pochhammer local system   of multiplicity one and rank $m=k\ell$.  
By Corollary \ref{gerad_proj_corol} we have
\[ \alpha_\infty =(g_{0,1}\,\cdots\,g_{0,\ell})\,( g_{1,1}\,\cdots\, g_{1,\ell})\,\cdots\,(g_{k-1,1}\,\cdots\, g_{k-1,\ell}),\] 
with $g_{j,i}=\omega_0^{k-1-j}\alpha_i\omega_0^{1-k+j}$, $j=0,\ldots,k-1$ and $i=1,\ldots,\ell$. 

Set $V_{j,i}=\mbox{ker}(\varrho_{\mathcal{K}}(g_{j,i})-\mathbf{I}_m)$. 
 Since  $\dim\ker(\varrho_{\mathcal{K}}(\alpha_\infty)-\mathbf{I}_m)=0$, 
 \begin{equation}\label{inters}  \bigcap_{i,j}V_{j,i}\subset \ker(\varrho_{\mathcal{K}}(\alpha_\infty)-\mathbf{I}_m)=(0). \end{equation}
Since $\mbox{codim}\ker(\varrho_{\mathcal{K}}(\alpha_i)-\mathbf{I}_m)=1$, $i=1,\ldots,\ell$, we get,  
by  similar arguments to Lemma  \ref{lemmag}  a unique decomposition into linear spaces of rank one, 
\[ \C^m=\bigoplus_{j=0}^{k-1}\bigoplus_{i=1}^\ell U_{j,i},\qquad V_{j,i}=\bigoplus_{(q,p)\neq (j,i)}U_{q,p}.  \] 
Now, we have the equivalences,  
\begin{eqnarray*}
u\in V_{j,i}&\Leftrightarrow& \varrho_{\mathcal{K}}\left(\omega_0^{k-j-1}\alpha_i\omega_0^{j-k+1}\right)(u)=u\\
            &\Leftrightarrow& \varrho_{\mathcal{K}}\left(\omega_0^{k-j}\alpha_i\omega_0^{j-k}\right)(\varrho_{\mathcal{K}}(\omega_0)(u))=\varrho_{\mathcal{K}}(\omega_0)(u)\\
            &\Leftrightarrow& \varrho_{\mathcal{K}}(\omega_0)(u)\in V_{j-1,i}.
\end{eqnarray*}
Hence $\varrho_{\mathcal{K}}(\omega_0)$ maps $V_{j,i}$ isomorphically onto $V_{j-1,i}$, for  $j=1,\ldots,k-1$ and $i=1,\ldots,\ell$.
Since  $\omega_0^k=\delta^r$, $\varrho_{\mathcal{K}}(\omega_0)$ also maps 
$V_{0,i}$  isomorphically onto $V_{k-1,i}$ for  $i=1,\ldots,\ell$.
Therefore, there is a basis $(u_{j,i})_{j,i}$  of $V$ such that   $U_{j,i}$ is the linear span of $u_{j,i}$ and the  matrix of $\varrho_{\mathcal{K}}(\omega_0)$ w.r.t.  $(u_{j,i})_{j,i}$
equals, 
\[     \left(\begin{array}{cccc} 
                                        & D_1\\
                                        &&\ddots\\
                                        &&&D_{k-1}\\
                                        D_0
				  \end{array}\right),
\]
with $D_0,\ldots,D_{k-1}$ diagonal matrices of order $\ell$, verifying $D_0\cdots D_{k-1}=\varepsilon^r\,\mathbf{I}_\ell$.
Up to a diagonal change of basis we can assume   $D_1=\ldots=D_{k-1}=I_\ell$ and $D_0=\varepsilon^rI_\ell$.
\quad$\Box$\vspace{2ex}

The proof of the previous lemma  showed  in particular,  that the monodromy along $\omega_0$ of an irreducible Pochhammer 
local system on a turbine without shaft has a very special shape: it  is always conjugated to $\oplus_i\zeta_i\mathbf{I}_\ell$ where   
 $\zeta_i$  ranges   through the set of   $k$-roots of $\varepsilon^r$ and $\mathbf{I}_\ell$ denotes the identity matrix of order $\ell$
(cf.  Theorem \ref{param_irred_without}). 

\subsection{Pochhammer local systems on the complement of an hypersurface}

 \vspace{1ex}

Let $Y$ be the germ of an hypersurface of a complex manifold $X$ and denote by   mult$(Y)$ its multiplicity. 
Let $Z$ be an irreducible component of $Y$.
Let $b\in Z$ be a nonsingular point of $Y$. 
Let $C$ be the germ at $b$ of a smooth curve transversal to $Z$.
Let $\gamma$ be a loop of $C\setminus Z$ with base point $c$.
Assume that $\oint_\gamma df/f=2\pi\sqrt{-1}$, where $f$ is a defining function of $Z$ at $b$.
Let $\mathcal L$ be a local system on $X\setminus Y$ and let  $M$ be its monodromy along the loop $\gamma$.
The conjugacy class $M_Z$ of $M$ does not depend on $\gamma$, $c$ or $b$ and is called the  \em local monodromy of $\mathcal L$ around $Z$. \em

Let $\mbox{rank}(\mathcal L)$ denote the dimension of the fiber of $\mathcal L$ at a point of $X\setminus Y$.
Let $j:X\setminus Y\hookrightarrow X$ be the inclusion map. 
The nonnegative  integer $\mbox{rank}(\mathcal L)-\dim(j_*\mathcal L)_b$ does not depend on $b$. 
It is called the \em multiplicity \em of $\mathcal L$ along $Z$ and it is denoted by mult$_Z(\mathcal L)$.
Notice that
\begin{equation}\label{MONOD}
\hbox{\rm mult}_Z(\mathcal L)=\hbox{\rm codim}(\ker(M_Z-id)).
\end{equation}

\begin{lemma}\label{gauss2}
Let $(Y,o)$ be the germ of an hypersurface in $\mathbb{C}^n$,  with irreducible components $Y_1,\ldots,Y_r$.
Let $X$ be an open neighbourhood of $o$.
Let $\tau: \mathbb C^n\to \mathbb C$ 
be a linear projection with fibers transversal to the tangent cone of $Y$.  
Let $\mathcal{L}$ be a local system on $X\setminus Y$.
Let $j$ be the open inclusion of $X\setminus Y$ into $X$.
The following statements are equivalent:
\begin{enumerate}
\item[$(a)$] $(j_*\mathcal L)_o=0$ and  $\hbox{\rm rank}(\mathcal L)=\sum_{i=1}^r\hbox{\rm mult}(Y_i)\hbox{\rm mult}_{Y_i}(\mathcal L)$; 
\item[$(b)$] For each $b\in \tau (U)$ such that $ b$ does not belong to the discriminant of $\tau$, $\mathcal L|_{\tau^{-1}(b)}$ 
verifies the conditions of  Lemma \ref{lemmag};
\item[$(c)$] $\R\tau_*(\mathcal L)$ vanishes.
\end{enumerate}
\end{lemma}
\noindent{\bf Proof.}
The equivalence between $(a)$ and $(b)$ follows from Lemma \ref{lemmag}.
The equivalence between $(b)$ and $(c)$ was proved in \cite{N}.\quad$\Box$\vspace{2ex}

\begin{definition} \cite{N} 
Let $(Y,o)$ be the germ of an hypersurface. 
Let $X$ be an open neighbourhood of $o$.
Let $\mathcal L$ be a local system on $X\setminus Y$.
We say that the local system $\mathcal L$ is  \em Pochhammer  \em if the conditions of Lemma \ref{gauss2} are verified.
\end{definition}







\section{Recognition of an irreducible   Pochhammer local system}\label{NR_turbine}

\subsection{An irreducibility  criterion for groups generated by pseudo-reflections}

\vspace{1ex}

Let   $\mathcal{H}$ be  a   linear subgroup of $\mbox{GL}_{m}(\C)$ 
generated by  pseudo-reflections  $A_i=\mathbf{I}_{m}-u_i\,v_i^T$, $i=1,\ldots,m$,  with  $u_i,v_i\in\C^m$, $u_i,v_i\neq 0$.

Let $\Gamma=(V,A)$ be a  digraph with set of vertices $V=\{x_1,\ldots,x_m\}$, and set of  arcs $A$
such that there is an arc from $x_i$ to $x_j$ ($i\neq j$),  i.e.,   $(x_i,x_j)\in A$,   if $v_i^Tu_j\neq 0$.

Recall that a digraph   is called  {\em strongly connected} if  for any distinct  pair of  vertices $x$ and $x'$ 
there is a directed path from  $x$ to $x'$. 

The following result gives an useful criterion to decide the  irreducibility  of   linear  subgroups     of   $\mbox{GL}_{m}(\C)$
generated by pseudo-reflections. 
\begin{theorem}{\em \cite[Theorem 5]{Fo}} \label{Fo}
The following hold.
\begin{enumerate}
\item[$(i)$] The group $\mathcal{H}$  is an irreducible subgroup   of $\mbox{GL}_m(\C)$   if and only if the following conditions are verified:
\begin{enumerate}
\item[$(a)$] $\Gamma$ is {\em strongly connected}.
\item[$(b)$] The matrix $[v_i^Tu_j]$, $i,j=1,\ldots,{m}$,  is invertible.
\end{enumerate}
\item[$(ii)$]   If  $\mathcal{H}$ is an irreducible subgroup   of $\mbox{GL}_m(\C)$,  there is a matrix $T\in\mbox{GL}_m(\C)$ such that  
$TA_iT^{-1}= \mathbf{I}_{m}-f_i\,e_i^T$, $i=1,\ldots,m$,  with $f_i\in\C^m$, $f_i\neq 0$, and  $e_i=(0,\ldots,0,1,0,\ldots,0)$ the $i$-th standard basis vector of  $\C^m$.
\end{enumerate}
\end{theorem}

When the pseudo-reflections are written in the `standard form' described in  $(ii)$, 
we  can replace the invertibility condition $(b)$   by  an invertibility condition expressed 
in terms of the product  $A:=A_1\cdots A_m$. More precisely,  we have the following result. 

\begin{lemma} \label{pseudo-reflection}
Consider pseudo-reflections $A_i= \mathbf{I}_{m}-f_i\,e_i^T$, with $f_i\in\C^m$, $f_i\neq 0$, for $i=1,\ldots,k$.  
Then $[e_i^Tf_j]$, $i,j=1,\ldots,{m}$, is invertible if and only if    $\mathbf{I}_{m}-A$ is invertible.  
\end{lemma}
{\bf Proof.}
For each $k=1,\ldots,m$,  we can find  complex numbers $\beta_{i,k}$, $i=1,\ldots,k-1$,  such that
\[ (\mathbf{I}_{m}-A)e_k=f_k+\sum_{i=1}^{k-1}\beta_{i,k}f_i.\] 
Hence  $\mathbf{I}_{m}-A$ is invertible  if and only if $f_1$, \ldots, $f_{m}$ are linearly independent.\quad$\Box$
  \vspace{2ex}

 \subsection{Numerical local data associated to an irreducible Pochhammer local system}

\vspace{1ex}

\begin{theorem}\label{param_irred}
 Let   $\mathcal{L}^*$ be an irreducible Pochhammer local system of multiplicity one on $T^*={{T}}^*_{n,k,\ell}$  with monodromy representation
 $\varrho:\pi_1(T^*)\to\mbox{GL}_m(\mathbb{C})$, $m=k\ell+1$. Assume moreover that $\mathcal{L}^*$ has semi-simple monodromy at $\infty$, i.e., 
  $\varrho(\omega_\infty)$  is semi-simple.
   Then the   following holds.
 \begin{enumerate}
\item[$(i^*)$]  $\varrho(\alpha_i)$,  are pseudo-reflections with special eigenvalues $\lambda_i\neq 1$, $i=0,\ldots,\ell$.
  \item[$(ii^*)$]  $\varrho(\delta)=\varepsilon\, \mathbf{I}_m$ for some $\varepsilon\in\C^*$.
 \item[$(iii^*)$]  $\varrho(\omega_0)$  is conjugated to $\zeta_1\mathbf{I}_\ell\oplus\cdots\oplus\zeta_{{k}}\mathbf{I}_\ell\oplus b$, 
 with $\{\zeta_1,\ldots,\zeta_{k}\}$  the set of ${k}$-roots of $\varepsilon^r$ and $b^k\neq\varepsilon^r$.
 \item[$(iv^*)$]  $\varrho(\omega_\infty)$ is conjugated to $\xi_1\mathbf{I}_{m_1}\oplus\cdots\oplus\xi_\nu\mathbf{I}_{m_\nu}$, with 
 $\xi_1,\ldots,\xi_\nu$, pairwise distinct nonzero complex numbers, $\xi_i^k\neq\varepsilon^r$, 
 $m_\nu\leq \cdots \leq m_1\leq \ell$ and  $\sum_jm_j=m$.
\item[$(v^*)$]\label{FuchsRelationwith}  $\lambda_1\,\cdots\,\lambda_\ell\,\xi_1^{m_1}\,\cdots\,\xi_{\nu}^{m_\nu}=(-1)^{(\nu-1)\ell}b\,\varepsilon^{r\ell}$ 
and $\lambda_0b^k=\varepsilon^r$.
\end{enumerate}
Furthermore,  $\mathcal{L}^*$ is rigid if and only if $\nu=k+1$,   $m_1=\cdots= m_k=\ell$ and $m_{k+1}=1$. 
\end{theorem}

We refer  to $(\varepsilon;b;\lambda_1,\ldots,\lambda_\ell;\xi_1,\ldots,\xi_\nu;m_1,\ldots,m_\nu)$
as the {\em numerical local data} of  $\mathcal{L}^*$.

We refer to  the former  relation of $(v^*)$ as the   {\em  Fuchs  relation}.

\vspace{1ex}

\noindent {\bf Proof.} $(i^*)$ follows  by definition. 
 
 $(ii^*)$  follows  from the  facts that  $\delta$ belongs to the center of $\pi_1(T^*)$ and $\mathcal{L}^*$ is irreducible.

Let us prove $(iii^*)$.
Since  each $\varrho(\alpha_i)$ is a pseudo-reflection,
$\mbox{codim}\ker(\varrho(\alpha_i)-\mathbf{I}_m)=1$ and we obtain, by a simple  linear algebra argument,
\begin{equation}\label{ell}
 \mbox{codim}\: \bigcap_{i=1}^\ell\ker(\varrho(\alpha_i)-\mathbf{I}_m)\leq\sum_{i=1}^\ell  \mbox{codim}\ker(\varrho(\alpha_i)-\mathbf{I}_m)=\ell.
 \end{equation}
Since $\alpha_1,\ldots,\alpha_\ell,\omega_0$ and $\delta$ generate $\pi_1(T^*)$, all eigenvalues of 
$\varrho(\omega_0)$ have geometric multiplicity less than or equal to $\ell$. 
Otherwise $\varrho$ would leave invariant  a nontrivial linear subspace, contradicting its irreducibility.  
 Since $\alpha_0\omega_0^k=\delta^r$ and $\varrho(\alpha_0)$ is a pseudo-reflection
with special eigenvalue $\lambda_0\neq  1$, $\varrho(\omega_0^k)$ is conjugated to $b'\oplus(\varepsilon^r)^{\oplus k\ell}$ 
with $b'\neq \varepsilon^r$. 
Hence  $\varrho(\omega_0)$  is   conjugated to  $\zeta_1\mathbf{I}_\ell\oplus\cdots\oplus\zeta_{{k}}\mathbf{I}_\ell\oplus b$ where $\zeta_1,\ldots,\zeta_{k}$
are the distinct  ${k}$-roots of $\varepsilon^r$ and    $b\neq \zeta_i$ for $i=1,\ldots,k$. 

Let us prove $(iv^*)$. 
Since   $\alpha_1,\ldots,\alpha_\ell,\omega_\infty$ and $\delta$  generate $\pi_1(T^*)$,  
all eigenvalues of  
$\varrho(\omega_\infty)$  have geometric multiplicity less than or equal to $\ell$
by an argument similar  to the argument used in  $(ii)$. By Corollary \ref{matrix-infty}  and the hypothesis, 
$\varrho(\alpha_\infty)=\varrho(\delta^r\omega_\infty^{-k})$ does not have the eigenvalue one. 
Hence   $\xi^k_i\neq\varepsilon^r$ for all $i$. 

$(v^*)$ follows  from   relations $\det(\varrho(\alpha_1)\cdots\varrho(\alpha_l)\varrho(\omega_\infty))=\det(\varrho(\omega_0))$ and
$\alpha_0=\delta^r\omega_0^{-k}$.

Assume that $\mathcal{L}^*$ is rigid. 
Let  $z(A)$  be the dimension of the centralizer  of a linear  endomorphism $A$. 
Since $\varrho(\delta)=\varepsilon\mathbf{I}_m$,  the pair $w_i,\delta$ generate $\pi_1(\partial N_i)$, $i=0,\infty$,   and the pair $\alpha_i,\delta$ 
generate $\pi_1(\partial N_i)$, $i=1,\ldots,\ell$,  we obtain
$\dim \Gamma(\partial N_i, \mbox{End}(\mathcal{L}^*))=z(\varrho(\omega_i))$, $i=0,\infty$ and 
$\dim \Gamma(\partial N_i, \mbox{End}(\mathcal{L}^*))=z(\varrho(\alpha_i))$, $i=1,\ldots,\ell$.
Since $\mathcal{L}^*$ is rigid, we have by  Theorem \ref{index-rigidity},
\[ \mbox{rig}(\mathcal{L}^*\,,\,T^*) =z(\varrho(\omega_\infty))+z(\varrho(\omega_0))+\sum_{i=1}^\ell z(\varrho(\alpha_i))-\ell(k\ell+1)^2=2. \]
 Thus  $z(\varrho(\omega_\infty))=k\ell^2+1$. Take $s\in\mathbb{Z}$ such that 
 $\ell=m_1=\cdots=m_s> m_{s+1}\geq\cdots\geq m_\nu\geq 1$.
Since  $\varrho(\omega_\infty)$ is semi-simple, $\sum_{i=1}^\nu m_i=k\ell+1$ and $\sum_{i=1}^\nu m^2_i=k\ell^2+1$. Hence 
\[ \sum_{i=s+1}^\nu m_i=(\nu-s)\ell+1\quad\mbox{and}\quad \sum_{i=s+1}^\nu m_i^2=(\nu-s)\ell^2+1. \]
Therefore,
\[ \ell\sum_{i=s+1}^\nu m_i=\sum_{i=s+1}^\nu m_i^2+\ell-1,\] 
that is, 
\[ \sum_{i=s+1}^\nu m_i(\ell-m_i)=\ell-1.\]
Since $t(\ell-t)>\ell-1$ if $t\in]1,\ell-1[$, $s=\nu-1$ and $m_\nu=1$. Therefore
$s=k$, $\nu={k}+1$, $m_1=m_2=\cdots=m_{{k}}=\ell$ and  $m_{{k}+1}=1$.\quad$\Box$
\vspace{2ex}


\begin{theorem}\label{param_irred_without} 
 Let  $\mathcal{L}$   be an irreducible Pochhammer local system of multiplicity one  on the turbine without shaft $T=T_{n,k,\ell}$ with monodromy representation
 $\varrho:\pi_1(T)\to\mbox{GL}_m(\mathbb{C})$, $m=k\ell$. Assume moreover that  $\mathcal{L}$ has semi-simple monodromy at infinity, i.e., 
 $\varrho(\omega_\infty)$  is semi-simple. 
 The  following holds. 
 \begin{enumerate}
\item[$(i)$] For each $i=1,\ldots,\ell$,  $\varrho(\alpha_i)$  is a  pseudo-reflection with special eigenvalue $\lambda_i\neq 1$. 
 \item[$(ii)$]  $\varrho(\delta)=\varepsilon\, \mathbf{I}_m$ for some $\varepsilon\in\C^*$.
 \item[$(iii)$]  $\varrho(\omega_0)$  is conjugated to $\zeta_1\mathbf{I}_\ell\oplus\cdots\oplus\zeta_{{k}}\mathbf{I}_\ell$, 
 with $\{\zeta_1,\ldots,\zeta_{k}\}$ the set of  ${k}$-roots of $\varepsilon^r$.
  \item[$(iv)$]  $\varrho(\omega_\infty)$ is conjugated to $\xi_1\mathbf{I}_{m_1}\oplus\cdots\oplus\xi_\nu\mathbf{I}_{m_\nu}$, with 
 $\xi_1,\ldots,\xi_\nu$, pairwise distinct nonzero complex numbers, $\xi_i^k\neq\varepsilon^r$, $m_\nu\leq \cdots \leq m_1\leq \ell$ and  $\sum_jm_j=m$.
\item[$(v)$]\label{FuchsRelationwithout}  
$\lambda_1\,\cdots\,\lambda_\ell\,\xi_1^{m_1}\,\cdots\,\xi_{\nu}^{m_\nu}=(-1)^{(\nu-1)\ell}\,\varepsilon^{r\ell}$. 
\end{enumerate}   
Moreover, $\mathcal{L}$ is rigid if and only if  $\nu=k+1$,  $m_1=\cdots= m_{k-1}=\ell$, $m_k=\ell-1$ and $m_{k+1}=1$. 
\end{theorem}

We  refer to $(\varepsilon;\lambda_1,\ldots,\lambda_\ell;\xi_1,\ldots,\xi_\nu;m_1,\ldots,m_\nu)$,  
as the  {\em numerical local data} of  $\mathcal{L}$. 

We refer to    relation  $(v)$ as the   {\em  Fuchs  relation}. 

\vspace{1ex}

\noindent{\bf Proof.} It follows   from similar arguments to the proof of  Theorem \ref{param_irred}, 
taking into account  Lemma \ref{matrix-infty}.  \quad $\Box$\vspace{2ex}







\section{Reconstruction  of a rigid Pochhammer local system}\label{recrec}

In this  section we solve the inverse problem of reconstructing  an rigid   irreducible  Pochhammer local system of multiplicity one on  a turbine from
its  numerical local data, assuming some additional generic conditions on this data. 

\subsection{The case $\ell=1$} 

\vspace{1ex}

\begin{theorem}\label{Reg-rec-with-irred}
A tuple  of nonzero complex numbers,  
\[ (\varepsilon;b;\lambda_1;\xi_1,\ldots,\xi_{k+1};1,\ldots,1),  \]
  verifying  Fuchs relation 
$\lambda_1\xi_1\,\cdots\,\xi_{k+1}=(-1)^{k}b\,\varepsilon^{r}$ 
 with  $\lambda_1\neq 1$  and  $b^k,\xi_j^k\neq\varepsilon^r$  for all $j$, plus   
   the generic condition  $\xi_i\neq b$ for all $i$, is  the    numerical local data of  an irreducible  rigid Pochhammer  local system  of multiplicity one on  $T^*_{n,k,1}$
and semi-simple monodromy at $\infty$. 
\end{theorem}
{\bf Proof.} Set $T^*=T^*_{n,k,1}$. 
In virtue of  theorems \ref{gerad_proj} and   \ref{param_irred},  it is enough   to find matrices $A_1,\Omega_0\in\mbox{GL}_m(\mathbb{C})$, $m=k+1$,  verifying:
\begin{enumerate}
\item[$(i^*)$]\label{eq-A_i_irred} $A_1$  is a pseudo-reflection with special eigenvalue  $\lambda_1$;
\item[$(ii^*)$]\label{eq-Omega_0_irred} $\Omega_0$  has eigenvalues $b,\zeta_1,\ldots, \zeta_{k}$,
 with $\zeta_1,\ldots,\zeta_{k}$  the distinct ${k}$-roots of $\varepsilon^r$ and $b\neq\zeta_i$ for all $i$;
\item[$(iii^*)$]\label{eq-A_0_irred}  $A_0=\varepsilon^r\Omega_0^{-k}$  is a pseudo-reflection with special eigenvalue  $\lambda_0=\varepsilon^rb^{-k}$;
\item[$(iv^*)$]\label{eq-Omega_infty_irred} $\Omega_\infty:=A_1^{-1}\Omega_0$ has $k+1$ distinct  eigenvalues $\xi_1,\ldots, \xi_{k+1}$;
\item[$(v^*)$]\label{eq_irred_rigid}  The linear representation  $\varrho:\pi_1(T^*)\to\mbox{GL}_m(\C)$, defined  
by $\varrho(\delta)=\varepsilon\mathbf{I}_m$, $\varrho(\alpha_i)=A_i$, $i=0,1$,  and 
$\varrho(\omega_0)=\Omega_0$, is   irreducible and rigid.
\end{enumerate} 
$(i^*)-(iii^*)$ are immediately fulfilled    considering,
\begin{equation}\label{matrices_with}
\Omega_0=\left(\begin{array}{c|cccc} b\\ \hline
                                        && 1\\
                                        &&&\ddots\\
                                        &&&&1\\
                                        1&\varepsilon^r
				  \end{array}\right),\qquad  
				   A_1=\left(\begin{array}{c|cccc} 1 &&&&a_k\\\hline&1&&& a_{k-1} \\ &&\ddots&&\vdots\\&&&1&a_1\\&&&&\lambda_1\end{array}\right),
\end{equation}
with $a_i\in\mathbb{C}$, $i=1,\ldots,k-1$, arbitrary (to be determined later). 

Let us prove $(iv^*)$. 
The matrix  $\Omega_\infty=A_1^{-1}\Omega_0$ equals 
\begin{equation}
\left(\begin{array}{c|cccc} b-\lambda_1^{-1}a_k&-\lambda_1^{-1}a_k\varepsilon^r\\ \hline
                                  -\lambda_1^{-1}a_{k-1}&-\lambda_1^{-1}a_{k-1}\varepsilon^r&1\\
                                  \vdots&\vdots&&\ddots& \\
                                  -\lambda_1^{-1}a_{1}&-\lambda_1^{-1}a_{1}\varepsilon^r&&&1\\
                                  \lambda_1^{-1}&\lambda_1^{-1}\varepsilon^r&&&
            \end{array}\right).
\end{equation}
Since  $\det(z\mathbf{I}_m-\Omega_\infty)=(b-z) q(z)-\lambda_1^{-1}a_kz^k$,  with
\[ q(z)=z^k+\varepsilon^r\lambda_1^{-1}\left(a_{k-1}z^{k-1}+a_{k-2}z^{k-2}+\cdots+a_1z-1\right),\]
 $\Omega_\infty$ has eigenvalues $\xi_1,\ldots, \xi_{k+1}$ if and only if  the  following conditions hold:
 \begin{eqnarray*}
-b\,\varepsilon^r\lambda_1^{-1}&=& (-1)^{k+1}\prod_{i=1}^{k+1}\xi_i,\\
   \varepsilon^r\lambda^{-1}(a_1b+1)&=&\sum_{\#\{i_1,i_2,\ldots,i_k\}=k}(-1)^k\xi_{i_1}\xi_{i_2}\cdots\xi_{i_k},\\
   \varepsilon^r\lambda^{-1}(a_2b-a_1)&=&\sum_{\#\{i_1,i_2,\ldots,i_{k-1}\}=k-1}(-1)^{k-1}\xi_{i_1}\xi_{i_2}\cdots\xi_{i_{k-1}},\\
   &\vdots&\\
      \varepsilon^r\lambda^{-1}(a_{k-1}b-a_{k-2})&=&\sum_{i_1\neq  i_2}\xi_{i_1}\xi_{i_2},\\
      b- \varepsilon^r\lambda^{-1}(a_k+a_{k-1})&=&-\sum_{i=1}^\nu\xi_{i}.\\
\end{eqnarray*}
First relation  is   equivalent to Fuchs relation.  The other relations determine the $a_i$'s  by  recurrence. 

Let us prove  $(v^*)$. The only nontrivial part is to prove the irreducibility. 

Set $M_0=A_0$ and set $M_j=\Omega_0^{k-j}A_1\Omega^{j-k}$,  $j=1,\ldots,k$. 
Each $M_j$ is a pseudo-reflection of the form $\mathbf{I}_k-v_j\,e_j^T$, for a unique nonzero  vector $v_j\in\mathbb{C}^{k+1}$.
 In order to prove that $\varrho$ is irreducible it is enough to show 
 that the linear subgroup $\mathcal{H}^*$  of $\mbox{GL}_{k+1}(\mathbb{C})$ generated by  $M_0,\ldots,M_k$ is irreducible. 

Let $\Gamma=(V,A)$ be the digraph with set of vertices $V=\{x_0,\ldots,x_k\}$ and set of arcs 
$A=\{(x_i,x_j)\::\:e_i^Tv_j\neq 0~(i\neq j)\}$.

 Since $b\neq \xi_i$ for all $i$ (by hypothesis),  $ -\lambda_1^{-1}a_kb^k=\prod_{i=1}^{k+1}(b-\xi_i)\neq 0$.  
Hence $a_k\neq 0$. Since   $e^T_0v_j=-b^{k-j}a_k\neq 0$,    $(x_0,x_j)\in A$ for every $j\neq 0$. 
Since   $e^T_j v_0=b^{k+1-j}\neq 0$,  $(x_j,x_0)\in A$ for every $j\neq 0$. Thus,   $\Gamma$ is strongly connected. 
Moreover, by Corollary \ref{gerad_proj_corol} and relation  $\xi_i^k\neq\varepsilon^r$,  $i=1,\ldots,k+1$, 
$M_0M_1\cdots M_k=\varepsilon^r\Omega_\infty^{-k}$ 
does not have the eigenvalue one. Hence $\mathcal{H}^*$  is irreducible by  Theorem \ref{Fo} and Lemma \ref{pseudo-reflection}.\quad$\Box$\vspace{2ex}



\begin{theorem}\label{Reg-rec-without-irred}
A tuple of nonzero complex numbers   $(\varepsilon,\lambda_1;\xi_1,\ldots,\xi_{k};1,\ldots,1)$ 
 verfying  Fuchs relation  $\lambda_1\xi_1\,\cdots\,\xi_{k}=(-1)^{k-1}\varepsilon^{r}$  
  with  $\lambda_1\neq 1$  and  $\xi_j^k\neq\varepsilon^r$  for all $j$,  
plus the generic  condition $\sum_{i=1}^k\xi_i\neq 0$, 
is  the    numerical local data of  an irreducible   Pochhammer   system  of multiplicity one on  $T_{n,k,1}$
with semi-simple monodromy at $\infty$. 
\end{theorem}
{\bf Proof.}
Set $T=T_{n,k,1}$. We have to find matrices $\Omega_0$ and $A_1$ such that   the following hold:
\begin{enumerate}
\item[$(i)$]\label{pt1} The eigenvalues of $\Omega_0$ are the  $k$-roots of $\varepsilon^r$;
\item[$(ii)$]\label{pt2}  $A_1$ is a pseudo-reflection with special eigenvalue $\lambda_1$;
\item[$(iii)$]\label{pt3}  $\Omega_\infty:=A_1^{-1}\Omega_0$ has distinct eigenvalues $\xi_1,\ldots, \xi_{k}$;
\item[$(iv)$]\label{pt4} The linear representation $\varrho:\pi_1(T)\to\mbox{GL}_m(\mathbb{C})$, 
 defined  by $\varrho(\delta)=\varepsilon\mathbf{I}_m$, $\varrho(\alpha_1)=A_1$ and  $\varrho(\omega_0)=\Omega_0$ is  irreducible and rigid.
\end{enumerate}
The conditions  ($i$)-($iii$) are fulfilled   at once     considering,  
\begin{equation}\label{matrices_without}
\Omega_0=\left(\begin{array}{cccc} &1\\ &&\ddots\\&&&1\\\varepsilon^r\end{array}\right), \qquad 
A_1=\left(\begin{array}{cccc} 1&&& b_{k-1} \\ &\ddots&&\vdots\\&&1&b_1\\&&&\lambda_1\end{array}\right) ,  
\end{equation}
with $b_i=\varepsilon^{-r}\lambda_1 a_i$, $i=1,\ldots,k-1$, such that  $z^k+\sum_{i=0}^{k-1}a_iz^i=\prod_{i=1}^{k}(z-\xi_i)$.

Let us  prove $(iv)$. The only nontrivial part is to prove that 
$\varrho$ is irreducible. 

Set $M_j=\Omega_0^{k-j}A_1\Omega_0^{j-k}$, $j=1,\ldots,k$. Then  $M_j=\mathbf{I}_k-v_j\,e_j^T$, for a unique nonzero
vector $v_j\in\mathbb{C}^k$ and it is enough to prove that the linear subgroup $\mathcal{H}$ of  $\mbox{GL}_k(\mathbb{C})$,  generated by the  $M_j$'s is irreducible. 

Let $\Gamma=(V,A)$ be the  digraph with set of vertices $V=\{x_1,\ldots,x_k\}$ and set of arcs 
$A=\{(x_i,x_j):e_i^Tv_j\neq 0~(i\neq j)\}$.

Since  $\Omega_0$  is the product   of an invertible diagonal matrix by   the  permutation matrix associated to the cycle
$\sigma=(1\:\cdots\:k)$, $M_{j}=\Omega_0M_{\sigma(j)}\Omega_0^{-1}$. Therefore   $(x_i,x_j)\in A$ 
 if and only if  $(x_{\sigma(i)},x_{\sigma(j)})\in A$.   
Since $\sum_{j=1}^k\xi_j\neq 0$,  $e_1^Tv_k=-b_{k-1}\neq 0$ and  thus  $(1,k)\in A$. Therefore $\Gamma$ is strongly connected.
By Corollary  \ref{gerad_proj_corol}   and the assumption $\xi_j^k\neq\varepsilon^r$ for all $j$, $M_1\cdots M_k=\Omega_0^k(\Omega_0^{-1}A_1)^k=\varepsilon^r\Omega_\infty^{-k}$,  
does not have the eigenvalue one.  Hence $\mathcal{H}$ is irreducible in virtue of   Theorem \ref{Fo} and Lemma \ref{pseudo-reflection}. \quad$\Box$
\vspace{2ex}

 \subsection{The general case} 

\vspace{1ex}

 Consider a turbine with shaft $T^*_{n,k,\ell}=N_\infty\setminus\bigcup_{i=1}^\ell\mbox{int}(N_i)$ and
  a solid torus $N_\bullet$ such that $ N_i\subset N_\bullet\subset N_\infty\setminus N_0$ for all $i$.
Set $T^*_{n,k,\bullet}:=N_\infty\setminus \mbox{int}(N_0\cup N_\bullet)$ and  $N'_\bullet=N_\bullet\setminus \cup_{i=1}^\ell \mbox{int}(N_i)$. 

\begin{figure}[ht]
\begin{center}
\input{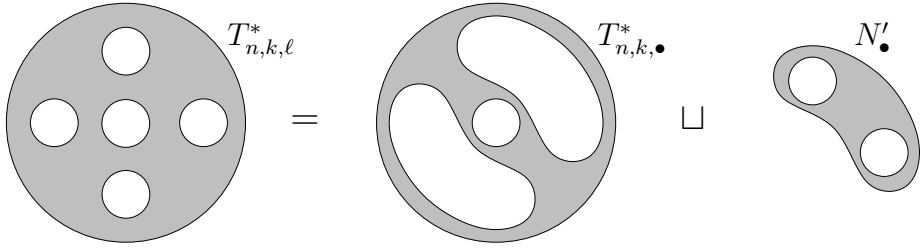}
\end{center}
\caption{\label{secturbines} Decomposition of a turbine with shaft}
\end{figure}

The topological space $T^*_{n,k,\bullet}$ defines a turbine with shaft and parameters $n,k,1$ and     $T^*_{n,k,\ell}$ is the union of   
$T^*_{n,k,\bullet}$ and  $N'_\bullet$,  pieced together along the torus $\partial  N_\bullet$ (see Figure \ref{secturbines}).

 Let $D_\bullet$ be a section transversal to $N_\bullet$. 
 Let $D_i\subset D_\bullet$ be a section transversal to $N_i$. 
Set  $D'_\bullet=D_\bullet\setminus\cup_{i=1}^\ell \mbox{int}(D_i)$. 
Set $U_\bullet=\partial D_\bullet$  and  $U_i=\partial D_i$, $i=1,\ldots,\ell$. Then 
$\pi_1(U_i)=\langle\alpha_i\rangle$ for all $i$,    $\pi_1(U_\bullet)=\langle\alpha_\bullet\rangle$  with 
 $\alpha_\bullet:=\alpha_1\cdots\alpha_\ell$, 
and  $\pi_1(D_\bullet^*)= \langle\alpha_\bullet,\alpha_1,\ldots,\alpha_\ell\::\:\alpha_\bullet=\alpha_1\cdots\alpha_\ell\rangle$. Moreover,
$\pi_1(N'_\bullet)\simeq\pi_1(D'_\bullet)\times \langle\delta\rangle$.

\begin{lemma}\label{decomp-pochh}
Let $\mathcal{L}$ be a  Pochhammer local system  on a turbine $T^*_{n,k,\ell}$ of multiplicity $\mu_i$ along $\partial N_i$, $i=0,\ldots,\ell$.  
The following holds:
\begin{enumerate}
\item $\mathcal{L}|_{T^*_{n,k,\bullet}}$ is a Pochhammer local system on $T^*_{n,k,\bullet}$  of multiplicity $k\mu_i$ along $\partial N_i$, $i=1,\ldots,\ell$, 
and multiplicity $\mu_0$ along $\partial N_0$.
\item $\mathcal{L}|_{D'_\bullet}\simeq\mathcal{K}\oplus\mathcal{\widetilde{L}}$, with 
 $\mathcal{K}$   a constant sheaf of rank $\mu_0+\sum_{i=1}^\ell(k-1)\mu_i$ and 
     $\widetilde{\mathcal{L}}$   a Pochhammer local system of multiplicity $\mu_i$ along $\alpha_i$.
\item The local data of $\mathcal{L}$ is determined by the local data of 
$\mathcal{L}'$ and of $\widetilde{\mathcal{L}}$. 
\end{enumerate}
\end{lemma}
\noindent{\bf Proof.} Straightforward.\quad$\Box$\vspace{2ex}

The previous lemma motivates the following definition.
\begin{definition}
Let $\mathcal{L}$, $\mathcal{L}'$  and  $\widetilde{\mathcal{L}}$   be  Pochhammer  
 local systems   on $T^*_{n,k,\ell}$,   $T^*_{n,k\bullet}$ and  $D'_\bullet$, respectively. 
We say that  $\mathcal{L}$   is an {\em extension}  of  $\mathcal{L}'$ via  $\widetilde{\mathcal{L}}$, and we denote it
by $\mathcal{L}=\mathcal{L}'\rtimes\widetilde{\mathcal{L}}$,   if  $\mathcal{L}|_{T^*_{n,k,\bullet}}\simeq \mathcal{L}'$ and there is constant sheaf $\mathcal{K}$  on $D'_\bullet$
such that  $\mathcal{L}|_{D'_\bullet}/\mathcal{K}\simeq \widetilde{\mathcal{L}}$.  
\end{definition}

All considerations and results above can be adapted, in a obvious way,   to   the turbines  without shaft, $T_{n,k,\ell}$ and  $T_{n,k,\bullet}$, 
where   $T_{n,k,\bullet}$  denotes the turbine with parameters $n,k,1$ obtained from $T^*_{n,k,\bullet}$ by `filling in' its shaft.


\begin{theorem}\label{Reg-rec-with}
Consider nonzero complex numbers $\varepsilon, b, \lambda_1,\ldots,\lambda_\ell, \xi_1,\ldots,\xi_{k+1}$,  such that 
 $\lambda_1,\ldots,\lambda_\ell,  \eta_1, \eta_2\neq 1$,  $\lambda_1,\ldots,\lambda_\ell\neq \eta_1$,  $\xi_i\neq \xi_j$ ($i\neq j$),
$\xi_1,\ldots,\xi_{k+1}\neq b$, $\xi_1^k,\,\ldots,\,\xi_{k+1}^k, \,b^k\neq\varepsilon^r$, $\sum_{i=1}^k\xi_i\neq 0$,  and $\eta_1\neq \eta_2$, where 
\[ \eta_1:=(-1)^{k-1}\varepsilon^r\prod_{i=1}^k \xi_i^{-1},\qquad  \eta_2:=(-1)^{k}\varepsilon^rb\prod_{i=1}^{k+1} \xi_i^{-1}. \]  
Assume moreover that the Fuchs relation  
\[  \lambda_1\,\cdots\,\lambda_\ell\,\xi_1^\ell\,\cdots\,\xi_{k}^\ell\,\xi_{k+1}=(-1)^{(k-1)\ell}b\,\varepsilon^{r\ell},
\]
holds. Then we have the following. 
\begin{enumerate}
\item[$(i)$]\label{F**} There is an  irreducible rigid Pochhammer local system  $\mathcal{F}$ of multiplicity one on $T_{n,k,\bullet}$,  that is 
determined by the  numerical local data $(\varepsilon;\eta_1;\xi_1,\ldots,\xi_k)$. 

\item[$(ii)$]\label{G^***} There is an irreducible rigid Pochhammer local system $\mathcal{G}^*$   of multiplicity one on $T^*_{n,k,\bullet}$,
that is determined by the 
numerical data $(\varepsilon;b;\eta_2;\xi_1,\ldots,\xi_{k+1})$.
\item[$(iii)$]\label{tildeX**} There is an irreducible  rigid   Pochhammer local   system $\widetilde{\mathcal{L}}$ of multiplicity one on $D'_\bullet$, 
that is determined by   the numerical  local data
 $(\lambda_1,\ldots, \lambda_\ell;\eta_1,\eta_2)$.
 \item[$(iv)$] The tuple $(\varepsilon;b;\lambda_1,\ldots,\lambda_\ell;\xi_1,\ldots,\xi_{k+1};\ell,\ldots,\ell,1)$ is  the    numerical local data of  an irreducible  rigid  
Pochhammer  local system  of multiplicity one on the turbine  $T^*_{n,k,\ell}$
and semi-simple monodromy at $\infty$,    obtained as extension of $\mathcal{F}^{\oplus(\ell-1)}|_{T^*_{n,k,\bullet}}\oplus\mathcal{G}^*$ via $\widetilde{\mathcal{L}}$. 
Moreover, this extension is unique. 
 \end{enumerate}
\end{theorem}
{\bf Proof.}  (i)-(iii) follow from  theorems \ref{Reg-rec-without-irred},  \ref{Reg-rec-with-irred} and  \ref{DSP2}.  Let us prove $(iv)$. 

Set $\mathcal{L}'=\mathcal{F}^{\oplus(\ell-1)}|_{T^*_{n,k,\bullet}}\oplus\mathcal{G}^*$ 
which is a Pochhammer local system on $T^*_\bullet$  of multiplicity $\ell$ along 
$\partial N_\bullet$ and multiplicity one along $\partial N_0$. Denote by $\varrho_{\mathcal{H}}$ the monodromy
representation of a local system $\mathcal{H}$.  By  (\ref{matrices_without}) and  (\ref{matrices_with}) we may assume that  
 $\varrho_{\mathcal{L}'}(\omega_0)=\varrho_{\mathcal{F}}(\omega_0)^{\oplus(\ell-1)}\oplus \varrho_{\mathcal{G}^*}(\omega_0)$ and
$\varrho_{\mathcal{L}'}(\alpha_\bullet)=\varrho_{\mathcal{F}}(\alpha_\bullet)^{\oplus(\ell-1)}\oplus \varrho_{\mathcal{G}^*}(\alpha_\bullet)$ 
are simultaneously conjugated, by a permutation matrix, to  the block  matrices 
\[ \Omega_0:=\left(\begin{array}{c|cccc}b\\ \hline& &\mathbf{I}_{\ell} \\ &&&\ddots\\&&&&\mathbf{I}_{\ell}\\  B&\varepsilon^r\mathbf{I}_\ell\end{array}\right),\qquad
A_\bullet:=\left(\begin{array}{c|cccc} 1 &&&&C_k\\\hline&\mathbf{I}_{\ell}&&& C_{k-1} \\ &&\ddots&&\vdots\\&&&\mathbf{I}_{\ell}&C_1\\&&&&C_0\end{array}\right),  \]
with $B=(0,\ldots,0,1)\in\C^\ell$, $C_0=\eta_1\,\mathbf{I}_{\ell-1}\oplus\eta_2$ , $C_i=b_i\,\mathbf{I}_{\ell-1}\oplus a_i$, $i=1,\ldots,k-1$,  
and $C^T_k=(0,\ldots,0,a_k)\in\C^\ell$, $a_k\neq 0$. 

By Theorem  \ref{DSP2} there is   an invertible matrix  $P$ of order $\ell$  such that $P^{-1}C_0P=A_1\cdots A_\ell$, where 
$A_i=\mathbf{I}_\ell-u_ie_i^T$ is a pseudo-reflection with special eigenvalue $\lambda_i$   and 
\begin{equation}\label{Ai} 
u_i=((\eta_1-\lambda_1)\eta_1^{-1},\ldots, (\eta_1-\lambda_{i-1})\eta_1^{-1},1-\lambda_i,\eta_1-\lambda_{i+1},\ldots,\eta_1-\lambda_{\ell}). 
\end{equation}
Set  $Q =1\oplus P^{\oplus k}\in\mbox{SL}_m(\mathbb{C})$ with $m=k\ell+1$. Then $\Omega'_0=Q^{-1}\Omega_0 Q$ and  $A'_\bullet=Q^{-1}A_\bullet Q$ 
are, respectively, given by
\[ \left(\begin{array}{c|cccc}b\\ \hline& &\mathbf{I}_{\ell} \\ &&&\ddots\\&&&&\mathbf{I}_{\ell}\\  B'&\varepsilon^r\mathbf{I}_\ell\end{array}\right),  \qquad \quad
   \left(\begin{array}{c|cccl} 1 &&&&C'_k\\\hline&\mathbf{I}_{\ell}&&& C'_{k-1} \\ &&\ddots&&\vdots\\&&&\mathbf{I}_{\ell}&C'_1\\&&&&C'_0\end{array}\right),  \] 
with $B'=P^{-1}B\neq 0$, $C'_j=P^{-1}C_jP$, $j=0,\ldots,k-1$ and   $C'_k=C_kP\neq 0$.

Set $\Theta=\{0\}\cup\bigcup_{j=0}^{k-1}\Theta_j$ with  $\Theta_j=\{(j,1),\ldots,(j,\ell)\}$.
Set $\nu(0)=0$ and  $\nu(j,i)=\ell j+i$ for each $(j,i)\in\Theta_j$. 
For each $\theta\in\Theta$, let $\hat{e}_\theta=(\hat{e}_{\theta,0},\ldots,\hat{e}_{\theta,k\ell})\in\mathbb{C}^m$, be the vector  
defined by  $\hat{e}_{\theta,r}=1$  if $r=\nu(\theta)$ and   $\hat{e}_{\theta,r}=0$, otherwise.   

Since $C'_0=A_1\ldots A_\ell$, there are unique    vectors  $v_{(k-1,i)}\in\mathbb{C}^m$, $i=1,\ldots,\ell$,  verifying 
\begin{equation}\label{vV}
\hat{e}^T_{(k-1,i)}v_{(k-1,i)}=e^T_iu_i,\qquad i=1,\ldots,\ell
\end{equation} 
such that   $A'_\bullet$ is the product of  pseudo-reflections  $X_i=\mathbf{I}_m-v_{(k-1,i)}\hat{e}_{(k-1,i)}^T$, $i=1,\ldots,\ell$. 
By Theorem \ref{gerad_proj}  we have  a linear representation, $\varrho:\pi_1({{T}}^*_{n,k,\ell})\to\mbox{GL}_m(\mathbb{C})$, 
defined by $\varrho(\delta)=\varepsilon\mathbf{I}_m$, $\varrho(\alpha_i)=X_i$, $i=1,\ldots,\ell$,  and   $\varrho(\omega_0)=\Omega'_0$.
By construction $\varrho$  is the monodromy representation  a Pochhammer local system    $\mathcal{L}^*$ extending $\mathcal{L}'$ via 
 $\widetilde{\mathcal{L}}$. 

It remains to  verify  that  $\mathcal{L}^*$   is irreducible and rigid.   The only nontrivial part is to prove the irreducibility.  
Set  $G_0=\varepsilon^r(\Omega'_0)^{-k}$ and $G_{(j,i)}=(\Omega'_0)^{k-1-j}X_i(\Omega'_0)^{1-k+j}$, $j=0,\ldots,k-1$.

Let $\Gamma=(V,A)$ be the digraph with set of vertices $V=\{x_\theta: \theta\in\Theta\}$ and set of arcs 
$A=\{(x_\theta,x_\eta):\hat{e}^T_\theta v_\eta\neq 0\}$. It is enough to prove that  $\Gamma$ is strongly connected.  
 
Let  $\Gamma_j$ be the induced sub-digraph of $\Gamma$, with set of vertices $V_j=\{x_\theta:\theta\in \Theta_j\}$.
 
By  (\ref{Ai}), (\ref{vV}),  and relations   $\lambda_i\neq\eta_1$, $i=1,\ldots,\ell$,  $\Gamma_{k-1}$ is the complete digraph. 
  
By   straightforward computations, reminding that $C'_k\neq 0$, 
 we have for   $j=0,\ldots,k-2$ and $i_1,i_2=1,\ldots,\ell$, 
 \[
\hat{e}^T_0\,v_{(j,1)}=b^{(k-j-1)}\,\hat{e}^T_0\,v_{(k-1,1)}\neq 0, \qquad \hat{e}^T_{(j,i_1)}\,v_{(j,i_2)}=\hat{e}^T_{(k-1,i_1)}\,v_{(k-1,i_2)}.  
\]
 Hence   $\Gamma_j\simeq\Gamma_{k-1}$   is   complete and    there are arcs  from vertex $x_0$ to  some    vertex of  $\Gamma_j$,  for every $j=0,\ldots,k-1$.
Since  $G_0=\varepsilon^rQ^{-1}\Omega_0^{-k}Q=I_m-v_0\,\hat{e}^T_0$ equals 
\[  \left(\begin{array}{c|ccc} \lambda_0 &&&\\\hline -B'b^{-k}&\mathbf{I}_{\ell}&& \\ \vdots&&\ddots&\\-B'b^{-1}&&&\mathbf{I}_{\ell}\end{array}\right),  \]
with $\lambda_0=\varepsilon^rb^{-k}$ and  $B'\neq 0$ there is $i_0\in\{1,\ldots,\ell\}$ such that $\hat{e}^T_{(j,i_0)}v_0\neq 0$, $j=0,\ldots,k-1$. 
Hence    there is  an  arc  from some  vertex of $\Gamma_j$ to     vertex $x_0$ for every $j=0,\ldots,k-1$. Therefore $\Gamma$ is strongly connected. 

By  Lemma \ref{decomp-pochh} every extension of $\mathcal{F}^{\oplus(\ell-1)}|_{T^*_{n,k,\bullet}}\oplus\mathcal{G}^*$ via $\widetilde{\mathcal{L}}$
have  local data isomorphic to the local data of $\mathcal{L}^*$. Hence these extensions are isomorphic to $\mathcal{L}^*$, since   $\mathcal{L}^*$ is  rigid.
\quad$\Box$\vspace{2ex}


\begin{theorem}\label{Reg-rec-without}
Consider nonzero complex numbers $\varepsilon,\lambda_1,\ldots,\lambda_\ell,\xi_1,\ldots,\xi_{k+1}$, with 
$\lambda_i,\eta_1,\eta_2$ distinct from one, $\lambda_1,\,\ldots,\,\lambda_\ell\neq \eta_1$,
$\xi_i\neq \xi_j$  $(i\neq j)$, $\xi_1+\cdots+\xi_k\neq 0$, $\xi_1+\cdots+\xi_{k-1}+\xi_{k+1}\neq 0$,  $\xi_1^k,\,\ldots,\,\xi_{k+1}^k\neq\varepsilon^r$, and  
where 
\[ \eta_1:=(-1)^{k-1}\varepsilon^r\prod_{i=1}^k \xi_i^{-1},\qquad  \eta_2:=(-1)^{k-1}\varepsilon^r\xi^{-1}_{k+1}\prod_{i=1}^{k-1} \xi_i^{-1}. \]
Assume moreover that the Fuchs relation   $\lambda_1\,\cdots\,\lambda_\ell\,\xi_1^\ell\,\cdots\,\xi_{k-1}^\ell\,\xi_{k}^{\ell-1}\,\xi_{k+1}=(-1)^{(k-1)\ell}\,\varepsilon^{r\ell}$ holds.
Then we have the following. 
\begin{enumerate}
\item[$(i)$]\label{F} There are   irreducible rigid Pochhammer local system   of multiplicity one on $T_{n,k,\bullet}$, $\mathcal{F}$  and  $\mathcal{G}$
determined, respectively,  by the 
numerical local data $(\varepsilon;\eta_1;\xi_1,\ldots,\xi_k)$ and  $(\varepsilon;\eta_2;\xi_1,\ldots,\xi_{k-1},\xi_{k+1})$.

\item[$(ii)$]\label{tildeX} There is a rigid   Pochhammer local   system $\widetilde{\mathcal{L}}$ of multiplicity one on $D'_\bullet$, determined by   the numerical  local data
$(\lambda_1,\ldots, \lambda_\ell;\eta_1,\eta_2)$.
\item[$(iii)$]  The tuple $(\varepsilon;\lambda_1,\ldots,\lambda_\ell;\xi_1,\ldots,\xi_{k+1};\ell,\ldots,\ell,\ell-1,1)$ is  the    numerical local data of  an irreducible  rigid  
Pochhammer  local system  of multiplicity one on the turbine  $T_{n,k,\ell}$
and semi-simple monodromy at infinity,    obtained as extension of $\mathcal{F}^{\oplus(\ell-1)}\oplus\mathcal{G}$ via $\widetilde{\mathcal{L}}$.
Moreover, this extension is unique. 
 \end{enumerate}
\end{theorem}
{\bf Proof.} 
 (i) and  (ii) follow from  theorems \ref{Reg-rec-without-irred} and  \ref{DSP2}.  Let us prove $(iii)$.

Set $\mathcal{L}'=\mathcal{F}^{\oplus(\ell-1)}\oplus\mathcal{G}$, which is a Pochhammer local system on ${{T}}_{n,k,\bullet}$ of multiplicity $\ell$.  Set 
\[ z^k+\sum_{i=0}^{k-1}a_iz^i=\prod_{i=1}^{k}(z-\xi_i),\quad z^k+\sum_{i=0}^{k-1}b_iz^i=(z-\xi_{k+1})\prod_{i=1}^{k-1}(z-\xi_i).\]
By  (\ref{matrices_without}) we can assume that $\varrho_{\mathcal{L}'}(\omega_0)$ and $\varrho_{\mathcal{L}'}(\alpha_\bullet)$ are simultaneously conjugated, by 
a permutation matrix, to  the block  matrices 
\[ \Omega_0=\left(\begin{array}{cccc} &\mathbf{I}_{\ell} \\ &&\ddots\\&&&\mathbf{I}_{\ell}\\\varepsilon^r\mathbf{I}_\ell\end{array}\right), \qquad 
A_\bullet=\left(\begin{array}{cccc} \mathbf{I}_{\ell}&&& C_{k-1} \\ &\ddots&&\vdots\\&&\mathbf{I}_{\ell}&C_1\\&&&C_0\end{array}\right) ,  \]
with $C_0=\eta_1\,\mathbf{I}_{\ell-1}\oplus\eta_2$ , $C_i=e_i\,\mathbf{I}_{\ell-1}\oplus f_i$, $i=1,\ldots,k-1$, 
where   $e_i=-a_0^{-1}a_i$ and $f_i=-b_0^{-1}b_i$, $i\geq 1$. Moreover, $C_{k-1}$ is in invertible matrix,  since
\begin{eqnarray*}
 e_{k-1}&=&(-1)^{k}\xi_1^{-1}\cdots\xi_k^{-1}(\xi_1+\cdots+\xi_k)\neq 0, \label{ek-1} \\
 f_{k-1}&=&(-1)^{k}\xi_1^{-1}\cdots\xi_{k-1}^{-1}\xi_{k+1}^{-1}(\xi_1+\cdots+\xi_{k-1}+\xi_{k+1})\neq 0 \label{fk-1}.
\end{eqnarray*}
Moreover,    $e_{k-1}\neq f_{k-1}$.
Let $A_i=\mathbf{I}_m-u_ie_i^T$, $i=1,\ldots,\ell$,  be  pseudo-reflections   with  $u_i$ given by (\ref{Ai}).  
Moreover,  $\widetilde{\mathcal{L}}$ is   determined by  the Pochhammer  tuple 
$(A_1,\ldots,A_\ell)$. 

By Theorem \ref{DSP2}, there is  an invertible  matrix $P\in\mbox{SL}_\ell(\mathbb{C})$ such that $P^{-1}C_0P=A_1\cdots A_\ell$. 
Set $Q=P^{\oplus k}\in\mbox{SL}_m(\mathbb{C})$.
Since $C_0^P=A_1\cdots A_\ell$,  $i=1,\ldots,\ell$, there are unique  vectors $v_{(k-1,i)}\in\mathbb{C}^m$, $i=1,\ldots,\ell$, 
such that  $A'_\bullet=Q^{-1}A_\bullet Q=X_1\cdots X_\ell$   with    $X_i=\mathbf{I}_m-v_{(k-1,i)}\hat{e}_{(k-1,i)}^T$ 
a  pseudo-reflection    with special 
eigenvalue $\lambda_i$. 
Furthermore, 
\begin{equation}\label{u-v-2}
\hat{e}^T_{(k-1,i)}v_{(k-1,i)}=e^T_iu_i,\qquad i=1,\ldots,\ell.
\end{equation} 
Hence we obtain a    linear representation  $\varrho:\pi_1({{T}}_{n,k,\ell})\to\mbox{GL}_m(\mathbb{C})$, 
defined by  $\varrho(\delta)=\varepsilon\mathbf{I}_m$, $\varrho(\alpha_i)=X_i$, $i=1,\ldots,\ell$,  and  $\varrho(\omega_0)=Q^{-1}\Omega_0Q=\Omega_0$. 
Moreover, by construction $\varrho$ is the
 monodromy representation of  a Pochhammer local system $\mathcal{L}$  of  multiplicity one extending $\mathcal{L}'$ to ${{T}}_{n,k,\ell}$  such that 
$\left(\mathcal{L}|_{U_\bullet}\right)/\mathcal{K}_{U_\bullet}=\widetilde{\mathcal{L}}$ with  $\mathcal{K}_{U_\bullet}$ 
a constant sheaf of rank $(k-1)\ell$. 
 
Let us prove that $\mathcal{L}$ is irreducible.

For $j=0,\ldots,k-1$ and $i=1,\ldots,\ell$,  set $G_{(j,i)}=(\Omega')_0^{k-1-j}X_i(\Omega'_0)^{1-k+j}$. 
Let $\Gamma=(V,A)$ be the digraph with set of vertices $V=\{x_{(j,i)}: 1\leq i\leq \ell,\: 0\leq j\leq k-1\}$ and set of arcs 
$A=\{(x_{(j,i)},x_{(\nu,\mu)}):\hat{e}^T_{(j,i)} v_{(\nu,\mu)}\neq 0\}$. It is enough to prove that the digraph $\Gamma$ is strongly connected.  
Let $\Gamma_j$ be the induced sub-digraph of $\Gamma$ 
with set of vertices $\{x_{(j,i)}:i=1,\ldots,\ell\}$. Since  entries of $u_i$ do not vanish, we have by (\ref{u-v-2}) that $\Gamma_{k-1}$ is  complete. 
Since  $\Omega_0$ is the product of an invertible   diagonal matrix by the permutation matrix associated to the 
 cycle $(1~\cdots~k)$, there are  nonzero complex numbers $\varepsilon_j$, $j=0,\ldots,k-1$, such that 
 \begin{equation}\label{theta}
 v_{(j,i)}=\varepsilon_j v_{(k-1,i)}, \qquad i=1,\ldots,\ell.
\end{equation}  
Hence  $\Gamma_j\simeq\Gamma_{k-1}$ are complete digraphs for all $j$. 
Since  $P^{-1}C_{k-1}P\neq 0$, there are  $i_1,i_2\in\{1,\ldots,\ell\}$ such that 
\[\hat{e}^T_{(0,i_1)}v_{(k-1,i_2)}\neq 0.  \]
Hence $((0,i_1),(k-1,i_2))\in A$. Therefore   $\Gamma$ is  strongly connected by  (\ref{theta}).

By   Lemma \ref{decomp-pochh}  every extension of $\mathcal{F}^{\oplus(\ell-1)}\oplus\mathcal{G}$ via $\widetilde{\mathcal{L}}$
have the same local data than $\mathcal{L}$. Hence the extension $\mathcal{L}$ is unique since it is rigid.
\quad$\Box$\vspace{2ex}

\begin{remark}
By Lemma \ref{equivalen} the irreducible rigid Pochhammer local systems on  $T^*_{n,k,\ell}$  constructed above  give rise to irreducible rigid local systems
on the Riemann sphere with  $\ell+2$ elements removed (which are no longer Pochhammer, except when $k=1$). The same conclusion holds replacing the turbine with shaft
$T^*_{n,k,\ell}$ by the turbine without shaft $T_{n,k,\ell}$. 
\end{remark}





\section{Application to the  weigthed homogeneous curves}\label{weighted_curves}\label{WH}

By  Theorem \ref{gerad_proj-curve},   the local fundamental group  of  $Y^*_{n,k,\ell}$ is isomorphic to the quotient group  of $\pi_1({{T}}^*_{n,k,\ell})$ 
by the relation $\omega_\infty^n=\delta^s$,  and henceforth  the monodromy at $\infty$ of    an irreducible local system on $Y^*_{n,k,\ell}$ is always semi-simple.  
By Theorem \ref{param_irred} the monodromy of an irreducible rigid Pochhammer local system of multiplicity one along the irreducible components of $Y^*$ is   
determined by a numerical tuple   $(\varepsilon;b;\lambda_1,\ldots,\lambda_\ell;\xi_1,\ldots,\xi_{k+1};\ell,\ldots,\ell,1)$  verifying  the assumptions  of  
Theorem  \ref{Reg-rec-with}, plus the additional condition 
\begin{equation}\label{curve}
\xi_i^n=\varepsilon^s,\qquad  \forall i.
\end{equation}
A similar remark apply to  the curve $Y_{n,k,\ell}$.

Keep the notations of  Section \ref{planecurves}. Choose  $\rho>0$ and $\ell+1$  nonzero complex numbers $C,C_1,\ldots,C_\ell$,  
such that  the  `fat curve'   $\widetilde{Y}_\bullet=\{(t^k,\lambda t^n): |\lambda^k-C|\leq\rho,\:t\in\C\}$,   contains 
the curves $Y_i$, $i=1,\ldots.\ell$. Set $Y_\bullet=\{y^k-Cx^n=0\}$ and $Y^*_\bullet=\{y(y^k-Cx^n)=0\}$.
\begin{theorem}\label{RigidLocalSystems_Y-with}
A tuple  $(\varepsilon;b;\lambda_1,\ldots,\lambda_\ell;\xi_1,\ldots,\xi_{k+1};\ell,\ldots,\ell,1)$  verifying  the assumptions  of  
Theorem  \ref{Reg-rec-with} and {\em (\ref{curve})},
is the numerical local data  of an  irreducible  rigid Pochhammer  local system $\mathcal{L}$  of multiplicity one along the   irreducible components of  $Y^*_{n,k,\ell}$. 
Moreover,  $\mathcal{L}|_{\mathbb{C}^2\setminus \widetilde{Y}^*_\bullet}\equiv\mathcal{F}^{\oplus\ell-1}|_{\mathbb{C}^2\setminus {Y}^*_\bullet}\oplus\mathcal{G}^*$,  
where  $\mathcal{F}$ is  the  rigid irreducible  Pocchammer local system of multiplicity one on $\mathbb{C}^2\setminus {Y}^*_\bullet$  determined by the tuple
 $(\varepsilon;\eta_1;\xi_1,\ldots,\xi_k;1,\ldots,1)$  and 
$\mathcal{G}^*$ is the rigid irreducible Pocchammer local system of multiplicity one on $\mathbb{C}^2\setminus {Y}_\bullet$ 
determined by  the tuple $(\varepsilon;b;\eta_2;\xi_1,\ldots,\xi_{k-1},\xi_k;1,\ldots,1)$.
\end{theorem}
\begin{theorem}\label{RigidLocalSystems_Y-without}
A tuple $(\varepsilon;\lambda_1,\ldots,\lambda_\ell;\xi_1,\ldots,\xi_{k+1};\ell,\ldots,\ell,\ell-1,1)$   verifying    the assumptions  of 
Theorem  \ref{Reg-rec-without} and  {\em (\ref{curve})},
is the numerical local data of an irreducible   rigid Pochhammer  local system of multiplicity one along the  irreducible components of  $Y_{n,k,\ell}$.
Moreover,  $\mathcal{L}|_{\mathbb{C}^2\setminus \widetilde{Y}_\bullet}\equiv\mathcal{F}^{\oplus\ell-1}\oplus\mathcal{G}$,  
with $\mathcal{F}$ and $\mathcal{G}$ rigid irreducible Pocchammer local systems of multiplicity one on $\mathbb{C}^2\setminus {Y}_\bullet$
 determined,  respectively, by the tuples
$(\varepsilon;\eta_1;\xi_1,\ldots,\xi_k;1,\ldots,1)$ and  $(\varepsilon;\eta_2;\xi_1,\ldots,\xi_{k-1},\xi_k;1,\ldots,1)$.
\end{theorem}


\begin{remark}
By  the main result of   \cite{LL},  there are vanishing sums of $k$ $n$-roots of the unity if  and only if $k\in \mathbb{Z}^+_0p_1+\cdots+\mathbb{Z}^+_0p_\nu$, where 
$p_1^{a_1}\cdots p_\nu^{a_\nu}$  is the prime factorization of $n$.  Therefore, when  $k\not\in\mathbb{Z}^+_0p_1+\cdots+\mathbb{Z}^+_0p_\nu$ the 
nonvanishing  assumptions on the sums of  $\xi_i$'s  can  be withdrawn  from 
theorems \ref{RigidLocalSystems_Y-with} and \ref{RigidLocalSystems_Y-without}.
 \end{remark}

\begin{remark}
 By  Theorem \ref{RigidLocalSystems_Y-without}, a tuple   $(\varepsilon;\lambda;\xi_1,\ldots,\xi_k)$, such that  
$\lambda\neq 1$, $\xi_i\neq\xi_j$  $(i\neq j)$, $\xi_i^n=\varepsilon^s$, $\sum_{i=1}^k\xi_i\neq 0$ and  $\lambda\,\xi_1\cdots\xi_k=(-1)^{k-1}\varepsilon^r$,
determines a rigid irreducible Pocchammer local system of multipliticity one along the cusp $Y_{n,k,1}=\{y^k=x^n\}$. 
Raising to the power $n$ the equality $\lambda\,\xi_1\cdots\xi_k=(-1)^{k-1}\varepsilon^r$ and reminding that $rn-sk=1$,  we get $\varepsilon=(-1)^{(k-1)n}\lambda^n$. 
 Thus we can replace the   tuple of data  above 
  by the   tuple of data $(\lambda;\zeta_1,\ldots,\zeta_k)$ where the $\zeta_i$'s are  $n$-roots of the unity verifying $\zeta_i\neq\zeta_j$ $(i\neq j)$,   $\zeta_1\cdots\zeta_k=1$ and $\sum_{i=1}^k\zeta_i\neq 0$ {\em(}cf. \cite{Si}{\em)}.
 \end{remark}

\subsection{An example of an irreducible  rigid Pochhammer  local system}\label{523}

\vspace{1ex}



Consider the germ of plane  curve at the origin, $(Y_{5,2,2},0)$, given by  $y^4-x^{10}=0$ and denote by $Y_i=\{y^2+(-1)^ix^5=0\}$, $i=1,2$,   
its   irreducible components, 
Consider a  tuple  of numerical local data 
$(\varepsilon\,;\,\lambda_1,\,\lambda_2\,;\,\xi_1,\,\xi_2,\,\xi_3;2,1,1)$ in   conditions of  Theorem  \ref{Reg-rec-without}, 
with $r=1$ and $s=2$.  Set  $\eta_1=-\varepsilon\xi_1^{-1}\xi_2^{-1}$ and $\eta_2=-\varepsilon\xi_1^{-1}\xi_3^{-1}$. 
By Fuchs relation, $\eta_1\eta_2=\lambda_1\lambda_2$.
Set $e_1=\xi_1^{-1}+\xi_2^{-1}$ and  $f_1=\xi_1^{-1}+\xi_3^{-1}$ and consider the matrices, 
  \[    \Omega_0=\left(\begin{array}{cc|cc} &&1&\phantom{x}  \\ &&&1 \\ \hline\varepsilon&\phantom{x}\\&\varepsilon &\phantom{x}\end{array}\right),\qquad\qquad
        A_*=\left(\begin{array}{cc|cc} 1&&e_1\\ &1&&f_1\\\hline&&\eta_1\\&&&\eta_2\end{array}\right).
\]
Let  
\[ A_1=\left(\begin{array}{cc} \lambda_1\\ \lambda_2-\eta_1 & 1\end{array}\right), \qquad\qquad 
A_2=\left(\begin{array}{cc}1& (\lambda_1-\eta_1)\eta^{-1}_1\\ &\lambda_2\end{array}\right), \] 
be the matrices given  by (\ref{A_i-Pochhammer}).  Let  
\[ P=\left(\begin{array}{cc} 1&1  \\ -\frac{\eta_1}{\lambda_1}& \frac{\lambda_2-\eta_1}{\lambda_1-\eta_1}\end{array}\right), \]
be the diagonalization matrix    of $A_1A_2$ associated to the  eigenvalues $\eta_1,\eta_2$. Set $Q=P^{\oplus 2}$. By straightforward computations we obtain, 
\[   QA_*Q^{-1} =  \left(\begin{array}{cccccc}
1&0&\frac{(\lambda_1-\eta_1)f_1-(\lambda_1-\eta_2)\eta_1}{\eta_2-\eta_1} &\frac{\lambda_1(\lambda_1-\eta_1)(f_1-e_1)}{\eta_1(\eta_2-\eta_1)} \\
0&1& \frac{(\lambda_2-\eta_1)(f_1-e_1)}{\eta_2-\eta_1}&\frac{\lambda_1(\eta_2-\lambda_1)f_1-(\eta_1-\lambda_1)e_1}{\eta_2-\eta_1}   \\
0&0&\lambda_1&\frac{\lambda_1(\lambda_1-\eta_1)}{\eta_1} \\
0&0& \lambda_2-\eta_1&-\lambda_1+\eta_2+\eta_1 
\end{array}\right).
\]
Replacing $\eta_1$, $\eta_2$, $e_1$ and $f_1$ by expressions only involving the numerical local data   we may write  the
 previous matrix as the product of  pseudo-reflections, 
\[   X_1=  \left(\begin{array}{cccc}
1&0&-\frac{\lambda_1-\varepsilon\xi^2_1}{\varepsilon\xi_1} &0\\ \\
0&1& -\frac{\lambda_2+\varepsilon\xi_1\xi_2}{\varepsilon\xi_1} &0 \\ \\
0&0&\lambda_1&0 \\ \\
0&0&\lambda_2+\varepsilon\xi_1\xi_2&1 
\end{array}\right),\qquad X_2=  \left(\begin{array}{cccc}
1&0&0&\frac{\lambda_1+\varepsilon\xi_1\xi_2}{\varepsilon\xi_2} \\ \\
0&1& 0&-\frac{\lambda_2-\varepsilon\xi_1^2}{\varepsilon\xi_1} \\ \\
0&0&1&-\frac{\lambda_1+\varepsilon\xi_1\xi_2}{\varepsilon\xi_1\xi_2} \\ \\
0&0& 0&\lambda_2 
\end{array}\right).
\]
The linear representation determined by   $\varrho(\alpha_i)=X_i$,  $i=1,2$, $\varrho(\omega_0)=\Omega_0$ and 
$\varrho(\delta)=\varepsilon\,\mathbf{I}_4$ is  the monodromy 
representation of an irreducible rigid  Pochhammer local system $\mathcal{L}$ of multiplicity one on the turbine  $T_{5,2,2}$ with the prescribed  
 local data 
$(\varepsilon;\lambda_1,\lambda_2;\xi_1,\xi_2,\xi_3;2,1,1)$. 
If we require additionally that  $\xi_i^5=\varepsilon^2$, $i=1,2,3$, then $\mathcal{L}$      determines an irreducible rigid   local system  
defined on complement of the weighted homogeneous curve $Y_{5,2,2}$. Furthermore, this local system is obtained as the  extension of a local system 
$\mathcal{L}_1\oplus\mathcal{L}_2$, where $\mathcal{L}_1$ and  $\mathcal{L}_2$ are irreducible rigid Pochhammer local systems on the complement of the cusp
$y^2=x^5$ determined by the tuples $(\varepsilon;\eta_1;\xi_1,\xi_2,1,1)$ and  $(\varepsilon;\eta_2;\xi_1,\xi_3;1,1)$, respectively.





\section{Special $\mathfrak D$-modules}\label{dmodules}

\subsection{Background on $\mathfrak{D}$-modules and perverse sheaves}

\vspace{1ex}

Let $X$ be a complex manifold of dimension $n$. 
Let $\mathcal O_X$ be the sheaf of holomorphic functions on  $T^*X$.
 Let $\pi:T^*X\to X$ be the \em cotangent bundle \em of $X$. 
Let $\theta$ be the canonical $1$-form of $X$. 
Let $(x_1,...,x_n)$ be a system of coordinates on an open set $U$ of $X$. 
Let us denote $x_i\circ \pi$ by $x_i$.
There is one and only one family of holomorphic functions $\xi_1,...,\xi_n$ on $\pi^{-1}(U)$ such that $\theta=\xi_1dx_1+\cdots+\xi_ndx_n$.
Moreover, $(x_1,...,x_n,\xi_1,...,\xi_n)$ is a system of coordinates on $\pi^{-1}(U)$.
We call {\em contact transformation} of $T^*X$ to a diffeomorphism between open sets of $T^*X$ that leaves $\theta$ invariant.
An analytic subset $\Gamma$ of $T^*X$ is called a (\em conic\em ) \em Lagrangean variety \em
if $\Gamma$ has dimension $n$ and the restriction of $\theta$ to the regular locus of $\Gamma$ vanishes.

A germ of a Lagrangean variety $\Gamma$ at a point $\alpha$ is said to be in  \em generic position \em if
$\Gamma\cap\pi^{-1}(\pi(\alpha))=\mathbb C\alpha$.

Let $Y$ be a closed analytic subset of $X$. 
The conormal of $Y$ is the smallest Lagrangean subvariety $\Gamma$ of $T^*X$ such that
$\pi(\Gamma)=Y$. 
We will denote $\Gamma$ by $T^*_YX$.
If $Y=\{ x\in U: ~x_1=\cdots=x_k=0\}$, $T^*_YX=\{ (x,\xi)\in \pi^{-1}(U):~ x_1=\cdots=x_k=\xi_{k+1}=\cdots =\xi_n=0\}$.
We identify $X$ with the {\em zero section} $T^*_XX$ of $T^*X$. 
If $Y$ is an hypersurface of $X$, $T^*_YX$ is 
the closure of the set of points $(a,\lambda df(a))$ such that $f(a)=0$, $df(a)\not=0$, 
$\lambda\in \mathbb C$ and $f$ is a local generator of the defining  ideal of $Y$.
We say that the germ of an hypersurface is in \em generic position \em if its conormal is in generic position.

Let $\mathbb C_X$ be the sheaf of locally constant functions on $X$.
A $\mathbb C_X$-module $F$ is called \em constructible \em if there is a decreasing sequence of closed analytic sets $(X_j)_{j\ge 0}$ 
such that $\cap_{j\ge 0}X_j=\emptyset$ and the sheaf $F|_{X_j\setminus X_{j+1}}$ is a local system for each $j$. 
A sheaf complex $F^\bullet:~\cdots \to F_k\to F_{k+1}\to \cdots$ is called \em constructible \em if its cohomology sheaves $H^j(F^\bullet)$ 
are constructible for each $j$ and $H^j(F^\bullet)=0$ for all but a finite number of $j$'s. We say that $F^\bullet$ verifies the \em conditions of support \em if 
 codim(supp$H^j(F^\bullet))\ge j$, for all $j$. 
 Let $\mathbb D^b_c(X)$ be the full subcategory of the derived category of sheaves $\mathbb D(X)$ whose objects 
 are the constructible complexes. Given a local system $G$ on $X$ the
  local system $G^{\curlyvee}=\mathcal H$om$_{\mathbb C_X}(G,\mathbb C_X)$ is called the \em dual \em of $G$.
The \em Poincar\'e-Verdier duality \em is a contravariant functor $\mathbb D_X:\mathbb D^b_c(X)\to \mathbb D^b_c(X)$ such that $\mathbb D_X(G)=G^\curlyvee$ for each local system $G$ on $X$. 
We say that an object $F$ of $\mathbb D^b_c(X)$ is a \em perverse sheaf \em if $F$ and its Poincar\'e-Verdier dual verify the conditions of support.
We will denote by $\mathbb Perv(X)$ the full subcategory of $\mathbb D^b_c(X)$ whose objects are perverse sheaves.

Let $\mathfrak D_X$ denote the ring of differential operators on $X$. Let $U$ be an open set of $X$.
Let $P=\sum_{|\alpha|\le m}a_\alpha\partial_x^\alpha$ be a differential operator of order $m$ over $U$.
Its principal symbol $\sigma(P)(x,\xi )=\sum_{|\alpha|= m}a_\alpha (x)\xi^\alpha$ is an holomorphic function homogeneous of degree $m$ on $\pi^{-1}(U)$.

Let $P=\sum_{i=0}^ma_i\partial^i_x\in \mathfrak D_\mathbb C$ be an ordinary differential operator. Assume that $a_m\not\equiv 0$. 
The $\mathfrak D_\mathbb C$-module $\mathfrak D_\mathbb C/\mathfrak D_\mathbb CP$ is a coherent $\mathfrak D_C$-module.
Moreover,
$$
\mathbb R\mathcal Hom_{\mathfrak D_\mathbb C}(\mathfrak{M},\mathcal O_\mathbb C) ~ : ~ 0 ~ \to ~\mathcal O_C ~ 
\stackrel{P}{\rightarrow} \mathcal O_C ~\to ~0.
$$
Hence $\mathcal Hom_{\mathfrak D_\mathbb C}(\mathfrak{M}, \mathcal O_\mathbb C)=\mathcal Ext^0_{\mathfrak D_\mathbb C}(\mathfrak{M},\mathcal O_\mathbb C)=Ker (P)$ and $\mathcal Ext^1_{\mathfrak D_\mathbb C}(\mathfrak{M},\mathcal O_\mathbb C)= Coker P$. Therefore $\mathcal Ext^0$ gives the solutions of the differential equation $Pf=0$ and $\mathcal Ext^1$ represents the obstructions to the existence of solutions of the inhomogeneous differential equation $Pf=g$.
In general, we call  $\mathcal Sol_X(\/-\/)=\mathbb R\mathcal Hom_{\mathfrak D_X}(\/-\/,\mathcal O_X)$ the \em solution functor. \em

Let $\Omega^\bullet_X$ be the differential forms complex of $X$. 
We call $d\mathcal R(\mathfrak{M})=\Omega^\bullet_X\otimes^\mathbb L_{\mathcal O_X}\mathfrak{M}$ the \em de Rham complex \em of $\mathfrak{M}$.
Notice that $d\mathcal R(\mathcal O_X)=\Omega^\bullet_X   \simeq\mathbb C_X$.
Moreover, $\mathbb D_X(\mathcal Sol_X(\mathfrak{M}))=d\mathcal R(\mathfrak{M})$.

We associate to a coherent $\mathfrak D_X$-module $\mathfrak{M}$ an analytic subset $\mathcal Char(\mathfrak{M})$ of $T^*X$, 
designated by  \em characteristic variety \em of $\mathfrak{M}$. 
If $\mathfrak{M}=\mathfrak D_Xu$, 
$$
\mathcal Char(\mathfrak{M}) ~ = ~ \{ (x,\xi): ~ \sigma (P)(x,\xi)= 0~ \hbox{\rm for all } ~ P ~ \hbox{ \rm such that } Pu=0 \}.
$$
Moreover, $\mathcal Char(\mathfrak D_\mathbb C/\mathfrak D_\mathbb CP)=\{ (x,\xi)\in T^*\mathbb C~: ~ \xi=0~$ or $ a_m(x)=0\}$.

\begin{theorem}
Let $\mathfrak{M}$ be a coherent $\mathfrak D_X$-module.  
The following statements are equivalent: $($a$)$ $\mathfrak{M}$ is a coherent $\mathcal O_X$-module; 
$($b$)$ $\mathfrak{M}$ is a locally free $\mathcal O_X$-module;  $($c$)$ $\mathcal Char(\mathfrak{M})=T^*_XX$;
$($d$)$ $\mathcal Sol_X(\mathfrak{M})$ is a local system; $($e$)$ $d\mathcal R(\mathfrak{M})$ is a local system.
\end{theorem}
A coherent $\mathfrak D_X$-module  $\mathfrak{M}$ verifying  the conditions above is 
called an \em integrable connection. \em
A coherent $\mathfrak D_X$-module $\mathfrak{M}$ is called \em holonomic \em if its characteristic variety is Lagrangean.
An integrable connection is the simplest example of a holonomic $\mathfrak D_X$-module.
If $\mathfrak{M}$ is a holonomic $\mathfrak D_X$-module, the hypersurface $Y=\pi(\mathcal Char(\mathfrak{M})\setminus X)$
is called the \em ramification locus \em of $\mathfrak{M}$. 
The singularities of the holomorphic solutions of $\mathfrak{M}$
are contained in its ramification locus. 

Let $\widehat{\mathcal O}_{X,a}$ denote the ring of formal power series of $X$ at the point $a\in X$.
A holonomic $\mathfrak D_X$-module is called \em regular holonomic \em if 
$$\mathbb R\mathcal Hom_{\mathfrak D_{X,a}}(\mathfrak{M}_a,\widehat{\mathcal O}_{X,a}/\mathcal O_{X,a})=0$$
for each $a\in X$. We will denote by ${\mathbb RH}(\mathfrak D_X)$ the category of regular holonomic $\mathfrak D_X$-modules.

The following theorem shows that regular holonomic $\mathfrak D$-modules are topological objects.

\begin{theorem}
The functor $\mathcal Sol$ \em [$d\mathcal R$] \em is a contravariant \em [\em covariant\em] \em equivalence of categories from 
${\mathbb RH}(\mathfrak D_X)$ onto $\mathbb Perv(X)$.
\end{theorem}

These equivalences of categories are known by the name of \em Riemann-Hilbert correspondence. \em 

We can associate to each constructible complex on $X$ an analytic subset of $T^*X$, the \em microssuport \em $SS(F)$ of $F$.
If $\mathfrak M$ is an holonomic $\mathfrak D_X$-module, 
$\mathcal Char (\mathfrak M)= SS(\mathcal Sol(\mathfrak M))=SS(d\mathcal R(\mathfrak M))$ (See \cite{KS}).

Let $\mathfrak{M}$ be a holonomic $\mathfrak D_X$-module.
We associate to each irreducible Lagrangean variety of $T^*X$  a nonnegative integer 
mult$_\Gamma(\mathfrak{M})$, the \em multiplicity of $\mathfrak{M}$ along $\Gamma$\em, 
such  that mult$_\Gamma(\mathfrak{M})\ge 1$ if and only if $\Gamma$ is an irreducible component of $\mathcal Char (\mathfrak{M})$.
If $0\to\mathfrak{M}'\to\mathfrak{M}\to\mathfrak{M}''\to 0$ is an exact sequence of $\mathfrak D_X$-modules 
and $\Gamma$ is an irreducible Lagrangean variety of $T^*X$, mult$_\Gamma (\mathfrak{M})=$mult$_\Gamma (\mathfrak{M}')+$mult$_\Gamma (\mathfrak{M}'')$.
If $\mathfrak{M}$ is an integrable connection on $X$, mult$_{T^*_XX}(\mathfrak{M})$ equals the rank of the $\mathcal O_X$-module $\mathfrak{M}$.
In particular,  mult$_{\{\xi=0\}}\mathfrak D/\mathfrak D P=m$. 
If $a_m(x)=x^l$, mult$_{\{x=0\}}$ $\mathfrak D/\mathfrak D P=l$.
Moreover,  Kashiwara's Index Theorem holds.
\begin{theorem} \em (See \cite{k2}, Theorem 6.3.1)\label{index}\em 
Let $\mathfrak{M}$ be a holonomic $\mathfrak D_X$-module. Let $(X_\al)$ be a Withney stractification of $X$ such that 
$\mathcal C$har$(\mathfrak{M})\subset \cup_\alpha T^*_{X_\alpha}X$. 
Let $d_\al$ be the codimension of $X_\al$ in $X$. 
For each $b\in X$ there are positive integers $C_b(\overline X_\al)$ such that
$$
\sum_i(-1)^i\dim H^i(\mathcal Sol(\mathfrak{M}))_b=  \sum_{b\in \overline X_\alpha} (-1)^{d_\al}C_b(\overline X_\al)
\mbox{mult}_{T^*_{X_\al}X}(\mathfrak{M}).
$$
If $Y$ is smooth at $b$, $C_b(Y)=1$. 
If $Y$ is a curve, $C_b(Y)$ equals the multiplicity of $Y$ at $b$.
\end{theorem}

In particular, if $\mathfrak M$ is an integrable connexion,
\begin{equation}\label{MMULTT}
\mbox{rank}\: \mathcal Sol (\mathfrak M) ~ = ~ 
\mbox{rank}\: H^0(d \mathcal R (\mathfrak M)) ~ = ~
\mbox{mult}_{T^*_XX}(\mathfrak M).
\end{equation}

Let $\mathfrak{M}$ be an holonomic $\mathfrak D_X$-module on a complex manifold $X$ of dimension $n$.
Let $\Omega_X^n$ denote the sheaf of differential forms of degree $n$ of $X$.
The holonomic $\mathfrak D_X$-module
\begin{equation}
\mathfrak{M}^* ~ = ~ \mathcal Ext_{\mathfrak D_X}^n(\mathfrak M, \mathfrak D_X)\otimes_{\mathcal O_X}\Omega^n
\end{equation}
is called
the \em dual \em of $\mathfrak{M}$. 
Notice that,
\[Ê\mathcal Char (\mathfrak{M}^*) = \mathcal Char (\mathfrak{M}),\qquad \mbox{mult}_\Gamma (\mathfrak{M})=mult_\Gamma(\mathfrak{M}^*),\]
 for each irreducible component of $\mathcal Char (\mathfrak{M})$ and
\begin{equation}\label{dual}
\mathcal Sol (\mathfrak{M}) ~ \simeq ~ d\mathcal R(\mathfrak{M}^*) ~ \simeq ~ \mathbb D (d\mathcal R(\mathfrak{M})).
\end{equation}

\subsection{Microlocal Riemann-Hilbert correspondence and special $\mathfrak{D}$-modules}


\vspace{1ex}

\begin{definition}\em 
Let $\mathfrak{M}$ be a germ of $\mathfrak D_X$-module  in generic position at a point $o$. 
Let $\alpha$  be the linear form such that $\pi^{-1}(o)\cap \mathcal Char(\mathfrak{M})=\mathbb C\alpha$.
We say that $\mathfrak{M}$ \em comes from an $\mathcal E_X$-module \em 
if there is a germ of a holomorphic vector field $\nu$ at $o$ such that $\sigma(\nu)(o)=\alpha$ and the map
$u\mapsto \nu u$  is a $\mathbb C$-linear isomorphism from $\mathfrak{M}_o$ onto $\mathfrak{M}_o$.
\end{definition}

\begin{definition}\em
Let $\mathfrak{M}$ be a $\mathfrak D_X$-module that comes from an $\mathcal E_X$-module.
We say that $\mathfrak{M}$ is a \em special $\mathfrak D_X$-module \em if 
$\mathfrak{M}^*$ also comes from an $\mathcal E_X$-module.
\end{definition}

For the motivation of the definition of $\mathfrak D_X$-module that comes from an $\mathcal E_X$-module see 
Theorem 8.6.19 of \cite{BJO}.

Let $X$ be a germ of a complex manifold at a point $o$.
Let $Y$ be a germ at $o$ a complex hypersurface with conormal in generic position.
Let $X_0$ be the complement in $X$ of the singular locus of $Y$.
Let $\imath: X_0\hookrightarrow X$ be the open inclusion.

We will denote by $\mathfrak D_\mu(X,Y)$ the category of germs at $o$ of regular holonomic $\mathfrak D_X$-modules
$\mathfrak M$ such that $\mathfrak M$ comes from an $\mathcal E_X$-module and the holomorphic solutions of $\mathfrak M$  ramify along $Y$.

We will denote by $\mathcal Special (X,Y)$ the subcategory of   $\mathfrak D_\mu(X,Y)$ consisting of the 
$\mathfrak D_X$-modules $\mathfrak M$ such that $\mathfrak M^*\in \mathfrak D_\mu(X,Y)$.

\begin{theorem}\label{RH}{\em \cite{N}}
The following statements hold:
\begin{enumerate}
\item[$(a)$]If $\mathfrak M$ is a holonomic $\mathfrak D_X$-module such that $\mathcal Char(\mathfrak M)\subset T_Y^*X\cup T^*_XX$,
$$H^i(d\mathcal R(\mathfrak M)) = H^i(\mathcal Sol(\mathfrak M)) = 0 ~\mbox{if} ~i\not=0,1.$$
\item[$(b)$]  The functor $\mathfrak{M} \mapsto d\mathcal R (\mathfrak{M})$ is an equivalence of categories
from $\mathfrak D_\mu(X,Y)$ onto the category of perverse sheaves $F$ such that $SS(F)$ equals $T^*_YX\cup T^*_XX$,
$F$ is concentrated in degree $0$, $F_o=0$ and $F|_{X\setminus Y}$ is a local system.
\item[$(c)$] If $\mathfrak M$ is in $\mathfrak D_\mu(X,Y)$,
$$
\mathcal Sol(\mathfrak M)|_{X_0} ~ \simeq ~ 
\mathbb D_{X_0} (\imath_*(d\mathcal R(\mathfrak M)|_{X\setminus Y})).
$$
\item [$(d)$] \label{RHH} The functor 
\[Ê\mathfrak{M} \mapsto d\mathcal R (\mathfrak{M})|_{X\setminus Y}, \]
defines an equivalence of categories between $\mathfrak D_\mu(X,Y)$ and the category of Pochhammer local systems on $X\setminus Y$.
\end{enumerate}
\end{theorem}

\begin{definition}\em
The equivalence of categories $(d)$ is called the \em Microlocal Riemann-Hilbert correspondence.\em
\end{definition}

\begin{theorem}\label{mmult}
Let $X$ be a complex manifold.
Let $\mathfrak{M}$ be the germ at a point $o$ of a holonomic $\mathfrak D_X$-module.
Assume that  $\mathfrak{M}$  comes from an $\mathcal E_X$-module.
Let $Y_i$, $1\le i\le \ell$, be the irreducible components of the ramification locus $Y$ of $\mathfrak{M}$. 
Let $\mathcal L$ be the Pochhammer local system associated to $\mathfrak{M}$ by the microlocal Riemann-Hilbert correspondence.
Then $\mathcal L$
is a local system of rank 
$$
\sum_{i=1}^\ell \mbox{mult}(Y_i) \mbox{ mult}_{T^*_{Y_i}X}(\mathfrak{M}).
$$ 
\em 
Moreover, \em mult$_{Y_i}(\mathcal L)$ $=$ mult$_{T^*_{Y_i}X}(\mathfrak{M})$, for $i=1,\ldots, \ell$.
\end{theorem}
{\bf Proof.}
Let $(S,o)$ be a germ of a smooth complex surface of $X$ that is noncharacteristic relatively to $\mathfrak{M}$ (see \cite{k2} or \cite{SS}).
Let $1\le i\le \ell$. 
We can assume that the canonical $1$-form of $T^*X$ equals $\eta dy+\sum_{r=1}^n\xi_rdx_r$, $Y_i$ is defined by a Weierstrass polynomial
$y^k+\sum_{j=0}^{k-1}a_j(x)y^j$ and $S$ equals $\{x_2=\cdots=x_n=0\}$. 
Hence $Y_i\cap S$ is a union of irreducible plane curves $Y_{ij}$ such that
\begin{equation}\label{UM}
\hbox{\rm mult}(Y_i)=\sum_j\hbox{\rm mult}(Y_{ij}).
\end{equation}
By (\ref{MONOD}), 
\begin{equation}\label{DOIS}
\hbox{\rm mult}_{Y_{i}}(d\mathcal R(\mathfrak{M})|_{X\setminus Y})=\hbox{\rm mult}_{Y_{ij}}(d\mathcal R(\mathfrak{M})|_{S\setminus Y}).
\end{equation}
Set $\Gamma_i=T^*_{Y_i}X$ and $\Gamma_{ij}=T^*_{Y_{ij}}S$.
Let us show that
\begin{equation}\label{TRES}
\hbox{\rm mult}_{T^*_{Y_{i}}X}(\mathfrak{M})=\hbox{\rm mult}_{T^*_{Y_{ij}}S}(\mathfrak{M}|_{S}).
\end{equation}
Since the multiplicity of $\mathcal C$har$(\mathfrak{M})|_S$ is calculated at a generic point of $\Gamma_{ij}$, we can assume that 
$\mathcal C$har$(\mathfrak{M})=\Gamma_i\cup T^*_XX$, 
$\mathcal C$har$(\mathfrak{M}|_S)=\Gamma_{ij}\cup T^*_SS$, 
\[
\Gamma_i=\{y=\xi_1=\cdots=\xi_n=0\},\qquad \Gamma_{ij}=\{(x_1,y;\xi_1,\eta): y=\xi_1=0\}.
\]
By Theorem 3.4.2 c) of \cite{SS} its enough to show that 
\begin{equation}\label{magic}
\rho_*\varpi^*[\Gamma_i]=
[\Gamma_{ij}],
\end{equation}
 where
$\varpi:S\times_{X}T^*X\hookrightarrow T^*X$ and $\rho: S\times_XT^*X\to T^*S$ is defined by
\[
\rho(x_1,y;\xi_1,\ldots,\xi_n,\eta)=(x_1,y,\xi_1,\eta).
\]
By Apendix D of \cite{SS},  in order to prove (\ref{magic})  it is enough to show that
\begin{equation}\label{magic1}
\mathcal O_{S\times_X T^*X}\otimes^{\mathbb L}_{\varpi^{-1}O_{T^*X}}
\varpi^{-1}\mathcal O_{T^*X}/I_{\Gamma_i}\simeq 
\mathcal O_{S\times_XT^*X}/T_{\Gamma_i\cap S\times_XT^*X}
\end{equation}
and
\begin{equation}\label{magic2}
\mathbb R\rho_*\mathcal O_{S\times_XT^*X}/I_{\Gamma_i\cap S\times_XT^*X}\simeq  \mathcal O_{T^*S}/I_{\Gamma_{ij}}.
\end{equation}

Let $\mathcal K^*$ be the Koszul complex of $\mathcal O_{T^*X}$ with respect to $y_1,\xi_1,\ldots,\xi_n$.
Since $\mathcal O_{S\times_XT^*X}\otimes_{\varpi^{-1}\mathcal O_{T^*X}}\varpi^{-1}\mathcal K $ is the Koszul complex of 
$\mathcal O_{S\times_XT^*X}$ with respect to $y_1,\xi_1,\ldots,\xi_n$ and
\[
\varpi^{-1}(\Gamma_i)=\{ (x_1,y,\xi_1,\ldots,\xi_n,\eta): y=\xi_1=\cdots=\xi_n=0\},
\]
\[
\mathcal O_{S\times_XT^*X}\otimes^{\mathbb L}_{\varpi^{-1}\mathcal O_{T^*X}}\varpi^{-1}\mathcal O_{T^*X}/I_{\Gamma_i}
~ \simeq ~
\]
\[
~ \simeq ~
\mathcal O_{S\times_XT^*X}\otimes_{\varpi^{-1}\mathcal O_{T^*X}}\varpi^{-1}\mathcal K 
~\simeq ~
\mathcal O_{S\times_XT^*X}/I_{\Gamma_{ij}}.
\]
Hence (\ref{magic1}) holds.

The map $\rho$ induces an homeomorphism $\widehat{\rho}$ from $\Gamma_i\cap S\times_XT^*X$ onto $\Gamma_{ij}$.
We wil denote by $\imath$ the maps $\Gamma_i\hookrightarrow S\times_{X}T^*X$, $\Gamma_{ij}\hookrightarrow T^*S$.
If $\mathcal I^*$ is a flabby resolution of $\mathcal O_{\Gamma_i\cap S\times_XT^*X}$,
\begin{equation}
\mathbb R\rho_*\mathcal O_{S\times_XT^*X}/I_{\Gamma_i\cap S\times_XT^*X} ~ \simeq ~ 
\rho_*\imath_*\mathcal I^*  ~ \simeq ~ 
\imath_*{\rho_0}_* \mathcal I^*.
\end{equation}
Since $\imath_*{\rho_0}_* \mathcal I^*$ is a flabby resolution of $\mathcal O_{T^*S}/I_{\Gamma_{ij}}$, (\ref{magic2}) holds.

By (\ref{UM}), (\ref{DOIS}) and (\ref{TRES}), we can assume that $X$ is a surface.
By Lemma 2.2 of \cite{N}, $H^j(d\mathcal R(\mathfrak{M}))=0$ for each $j\not=0,1$. 
By statement (3) of \cite{N}, 
$H^j(d\mathcal R(\mathfrak{M}))_o=0$ for each $j$.
By (\ref{dual}) we can replace in Theorem \ref{index} $\mathcal Sol(\mathfrak{M})$ by $d\mathcal R(\mathfrak{M})$.
 Applying Kashiwara's Index Theorem at the points  
 $o$,  $b\in X\setminus Y$ and $b_i\in Y_i\setminus \{o\}$, $1\le i\le \ell$,
 we conclude that
\begin{equation}\label{igualdade1}
\dim H^0(d\mathcal R(\mathfrak{M}))_b~=
~\textrm{mult}_{T^*_XX}(\mathfrak{M})~=
~\sum_{i=1}^\ell\textrm{mult}_o(Y_i)~\textrm{mult}_{T^*_{Y_i}X}(\mathfrak{M}),
\end{equation}
and
\begin{equation}\label{igualdade2}
\textrm{mult}_{T^*_XX}(\mathfrak{M})-\textrm{mult}_{T^*_{Y_i}X}(\mathfrak{M})~=
~\dim H^0(d\mathcal R(\mathfrak{M}))_{b_i}-\dim H^1(d\mathcal R(\mathfrak{M}))_{b_i}.
\end{equation}
By Theorem \ref{RH}, $H^1(d\mathcal R(\mathfrak{M}))_{b_i}$ vanishes for  $1\le i\le \ell$.\quad$\Box$\vspace{2ex}

Let $\mathcal Pochh^*(X,Y)$ be the category of Pochhammer local systems on $X\setminus Y$ such that
the Verdier dual of $\imath_* d\mathcal R(\mathfrak M)|_{X\setminus Y}$ is concentrated in degree $0$.

\begin{theorem}\label{mourinhos}
Let $\mathfrak M$ be a $\mathfrak D_X$-module that comes from an $\mathcal E_X$-module.
The following statements are equivalent:
\begin{enumerate}
\item[(a)] 
$\mathfrak M$ is a special $\mathfrak D_X$-module.

\item[(b)]  $\mathcal Ext_{\mathfrak D_X}^1(\mathfrak M,\mathcal O)$ vanishes.

\item[(c)]  $d\mathcal R(\mathfrak M)|_{X\setminus Y}$ is in $\mathcal Pochh^*(X,Y)$.
\end{enumerate}
\end{theorem}
{\bf Proof.}
If $\mathfrak M$ is special, (b) follows from (\ref{dual}).

Let us show that $\mathcal M$ is special when (b) holds.
By Theorem \ref{RH} (b), it is enough to show that $\mathcal Hom_{\mathfrak D}(\mathfrak M,\mathcal O)_o$ vanishes.
By Kashiwara's Index Theorem, dim$\mathcal Hom_{\mathfrak D_X}(\mathfrak M,\mathcal O)_o$ equals
\begin{equation}\label{zero}
\mbox{mult}_{T^*_XX}(\mathfrak M)-\sum_{i=1}^\ell C_o(Y_i)\mbox{mult}_{T^*_{Y_i}X}(\mathfrak M).
\end{equation}
By the arguments of the proof of Theorem \ref{mmult}, we can assume that $Y$ is a plane curve. 
Since $\mathfrak M$ comes from an $\mathcal E_X$-module, it follows from Theorem \ref{mmult} that (\ref{zero}) vanishes. 

The equivalence between statements (b) and (c) follows from Theorem \ref{RH} (c).\quad$\Box$\vspace{2ex}

\begin{corollary}
The functors $\mathfrak M\mapsto d\mathcal R(\mathfrak M)|_{X\setminus Y}$ and
$\mathfrak M\mapsto \mathcal Sol(\mathfrak M)|_{X\setminus Y}$
define equivalences of categories between $\mathcal Special (X,Y)$ and $\mathcal Pochh^* (X,Y)$.
\end{corollary}

\begin{theorem}\label{mourinho}
Let $\mathfrak{M}$ be a $\mathfrak D_X$-module that comes from an $\mathcal E$-module, with ramification locus $Y$.
Let $Y_i$, $1\le i\le \ell$, be the irreducible componentes of $Y$.
Let $\mathcal L$ be the Pochhammer local system associated to $\mathfrak{M}$ by the Riemann-Hilbert correspondence.
Let $\sigma_i$ be the sum of the geometric multiplicities of the eigenvalues of the local monodromy of $\mathcal L$ around $Y_i$, $1\le i\le \ell$.

Then $\sigma_i\le$mult$_{Y_i}(\mathfrak{M})$, $1\le i\le \ell$.
If $\sigma_i=$mult$_{Y_i}(\mathfrak{M})$, $1\le i\le \ell$,  
$\mathfrak{M}$ is a special $\mathfrak D$-module.
\end{theorem}
{\bf Proof.}
Let $1\le i\le \ell$. Since $\sigma_i\le$mult$_{Y_i}(\mathcal L)$, Theorem \ref{mmult} implies  that $\sigma_i\le$mult$_{T^*_{Y_i}X}(\mathfrak{M})$.

Assume that $\sigma_i$ equals {mult}$_{T^*_{Y_i}X}(\mathfrak{M})$ for each $i$.
Applying the Kashiwara index theorem at a regular point $b_i$ of $Y$ that belongs to $Y_i$, we have, 
\[Ê\mbox{dim}H^0(\mathcal Sol(\mathfrak{M}))_{b_i} ~ = ~ 
\mbox{mult}_{T^*_XX}(\mathfrak{M})-\mbox{mult}_{T^*_{Y_i}X}(\mathfrak{M})+\mbox{dim}H^1(\mathcal Sol (\mathfrak{M}))_{b_i}. \]
Since 
\[ \mbox{dim}H^0(\mathcal S (\mathfrak{M}))_{b_i}\le H^0(\mathcal Sol(\mathfrak{M}))_b-\sigma_i=
\mbox{mult}_{T^*_XX}(\mathfrak{M})-\mbox{mult}_{T^*_{Y_i}X}(\mathfrak{M}),\]
$\mathcal Ext^1_\mathfrak D(\mathfrak{M},\mathcal O)$ vanishes outside the singular locus of $Y$. By the support conditions for perverse sheaves, 
$\mathcal Ext^1_\mathfrak D(\mathfrak{M},\mathcal O)$ vanishes.\quad$\Box$\vspace{2ex}

Let $\mathcal F$ be a $\mathfrak D_X$-module and let $\varphi$ be a section of $\mathcal F$.
Let $\mathcal J$ be the left ideal consisting of  sections $P$ of $\mathfrak D_X$ such that $P\varphi=0$.
We will identify the $\mathfrak D_X$-modules $\mathfrak{M}=\mathfrak D_X/\mathcal J$ and $\mathfrak D_X\varphi$.
If $\mathfrak{M}$ is a coherent $\mathfrak D_X$-module, we will say that $\mathfrak{M}$ is {\em generated by its solution} $\varphi$.

Let $\Omega$ be a polydisc of $\mathbb C^n$ centered at $0$.
Let $Y$ be a germ of hypersurface of $(\mathbb C^n,0)$.
Let $\tau_0:\widetilde\Omega\to \Omega\setminus Y$ be the universal covering of $\Omega\setminus Y$.
Let $j:\Omega \setminus Y\hookrightarrow \Omega$ be the open inclusion.
Set $\tau=j\tau_0$. 
The sheaf $\widetilde{\mathcal O}=\tau_*\mathcal O_{\widetilde\Omega}$ is the called the 
\em sheaf of multivalued holomorphic functions of $\Omega$ ramified along $Y$. \em 
The sheaves $\widetilde{\mathcal O}$ and $\widetilde{\mathcal O}/\mathcal O$ are $\mathfrak D_\Omega$-modules.

\begin{theorem}\label{speciald}
Let $Y$ be a germ of an hypersurface of a complex manifold $X$.
Let $\mathfrak{M}$ be a special $\mathfrak D_X$-module with ramification locus $Y$. 
There is a germ of a multivalued holomorphic function  $f$ ramified along $Y$ such that $\mathfrak{M}=\mathfrak D_Xf$.
\end{theorem}
{\bf Proof.}
By theorems 4.5.2 and 5.1.1 of \cite{KK}, there is $\varphi\in \widetilde{\mathcal O}/\mathcal O$ such that 
$\mathfrak{M}=\mathfrak D_X\varphi$. 
There are differential operators $P_1,\ldots, P_m$ such that $\mathfrak{M}=\mathfrak D_X/(\mathfrak D_XP_1+\cdots+\mathfrak D_XP_m)$.
Moreover, there is a free resolution
\begin{equation}\label{resolut1}
\mathfrak D_X^n   \rightarrow{(P_{ji})} \mathfrak D_X^m \rightarrow{(P_i)} \mathfrak{M} \to  0.
\end{equation}
Applying the functor $\mathcal Hom_{\mathfrak D_X}(*,\mathcal O)$ to (\ref{resolut1}) we obtain the exact sequence
\begin{equation}\label{resolut2}
0 \to \mathcal H om_{\mathfrak D_X}(\mathfrak{M},\mathcal O)  \rightarrow{(P_j)} \mathcal O^m \rightarrow{(P_{ij})} \mathcal O^n.  
\end{equation}

\noindent
Let $f_0\in\widetilde{\mathcal O}$ be a representative of $\varphi$.
Since $P_jf_0\in\mathcal O$, $1\le j\le m$, $\sum_{j=1}^mP_{ij}(P_jf_0)=0$, $i=1,\ldots,n$, and $\mathcal Ext^1_{\mathfrak D_X}(\mathfrak{M},\mathcal O)$ vanishes, 
there is $g\in\mathcal O$ such that $P_jg=P_jf_0$, $1\le j\le m$.
We take $f=f_0-g$.\quad$\Box$\vspace{2ex}

\begin{definition}\label{specialf}\em
We say that the germ of a multivalued holomorphic function is a \em special function \em if it generates a special $\mathfrak D_X$-module.
\end{definition}


\subsubsection{Examples}

Consider $X=\mathbb{C}$ and set $\mathfrak{M}=\mathfrak D_{\mathbb C}/\mathfrak D_{\mathbb C}x\partial_x$. 
The $\mathfrak D_{\mathbb C}$-module $\mathfrak{M}$ is generated by its holomorphic microfunction solution $\log x+\mathcal O_\mathbb C$. 
Remark that 
\begin{equation}
\mathcal Ext^1_{\mathfrak D_\mathbb C}(\mathfrak{M},\mathcal O_\mathbb C)= Coker ~~ x\partial_x:\mathcal O_\mathbb C \to \mathcal O_\mathbb C
\end{equation}
 does not vanish. 
The only multivalued holomorphic solutions of $\mathfrak{M}$ are the constant functions. 
Since the germ of $\mathfrak D_{\mathbb C}1$ at the origin equals $\mathbb C\{x\}$ and $\mathfrak{M}_0$ is isomorphic  to $\mathbb C\{x\}\oplus \partial_x\mathbb C[\partial_x]$, there is no germ of multivalued holomorphic solution of $\mathfrak{M}$ that generates $\mathfrak{M}$.


Consider now  $X=\mathbb C^m$ with coordinates $(x,y,t_1,\ldots,t_{m-2})$.
Set $Y=\{ (x,y,t) :y^k-x^n=0\}$, where $(k,n)=1$ and $2\le k\le n-1$.
Set $\vartheta=x\partial_x+(n/k)\partial_y$. Given complex numbers
$\lambda_i$, $i\in \mathbb Z$, such that
\[
\lambda_{i+k}=\lambda_i, \qquad\mbox{and}\qquad  \alpha_i:=\lambda_i-\lambda_{i+1}+\frac{n-k}{k}
\]
is a nonnegative integer, we will denote by
$\mathfrak{M}_{(\lambda_i)}$ the $\mathfrak D_X$-module given by 
generators $u_i$, $i\in \mathbb Z$, and relations
\[
u_{i+k}=u_i,\qquad (\vartheta-\lambda_i)u_i=0, \quad
\partial_xu_i=-\frac{n}{k}x^{\alpha_i}\partial_yu_{i+1}, \quad \partial_{t_j}u_i=0,~~ \forall j.
\]
The systems $\mathfrak M_{(\lambda_i)}$ have characteristic  variety
$T^*_YX\cup T^*_XX$ and multiplicity $1$ along $T^*_YX$, and henceforth are regular holonomic. 
By Theorem 5.2 of  \cite{NeSi}   there is a
denumerable subset $\Xi$ of $\mathbb C$ such that  $\mathfrak
M_{(\lambda_i)}$  is special if and only if
$\lambda_i\not\in\Xi$, $ i\in \mathbb Z$.

By Theorem 4.6 of \cite{NeSi} the $\mathfrak D_X$-modules $\mathfrak
M_{(\lambda_i)}$ verifying the conditions above are the only special
$\mathfrak D$-modules with characteristic variety $T^*_YX\cup T^*_XX$
and multiplicity $1$ along  $T^*_YX$.
By Theorem 5.3 of \cite{NeSi} $\mathfrak M_{(\lambda_i)}$ is generated
by a multivalued holomorphic function of the type
\[
y^a\left( \frac{y^k}{x^n}  \right)^bF\left(  \frac{y^k}{x^n}\right),
\]
where $F$ is a $_kF_{k-1}$ hypergeometric function.

\vspace{2ex}

Next result follows from theorems  \ref{RigidLocalSystems_Y-with}, \ref{RigidLocalSystems_Y-without}, \ref{RH} and \ref{mourinho}.
\begin{theorem}
Let $Y$ be a weigthed homogeneous plane curve with irreducible components $Y_i$, $i\in I$.
Let $\lambda_i$, $i\in I$, be complex numbers such that $\lambda_i\not=0,1$ .
There is an irreducible special $\mathfrak D$-module $\mathfrak{M}$ with solutions ramified along $Y$ such that the following hold  for each $i\in I$:
\begin{enumerate}
\item[(a)]  $\mathfrak{M}$ has multiplicity $1$ along the conormal of $Y_i$.
\item[(b)] The special eigenvalue of the local monodromy of $d\mathcal R(\mathfrak{M})$ around $Y_i$ equals $\lambda_i$.
\end{enumerate}
\end{theorem}

\vspace{2ex}

\section*{Acknowledgements}
 The authors thank Pierre Deligne   \cite{De} for his valuable remarks on \cite{Si}.

\end{document}